\documentclass[12pt]{article}

\title{\bf Equidistants  for families of surfaces}
\author{Peter Giblin and Graham Reeve }

\usepackage{amssymb}
\usepackage{amsmath}
\usepackage{amsfonts}
\usepackage{graphicx}
\usepackage{color}
\usepackage[all]{xy}
\usepackage[width=\textwidth]{caption}

\date{}
\newcommand{\RR}{\mathbb R}
\newcommand{\cS}{\mathcal S}

\newcommand{\g}{\gamma}

\newcommand{\cA}{\mathcal A}

\newcommand{\ee}{\varepsilon}

\newcommand{\wH}{\widetilde{H}}  
\newcommand{\wh}{\tilde{h}}

\newcommand{\gl}{\lambda}

\newcommand{\bp}{\mbox{\boldmath $p$}}
\newcommand{\bq}{\mbox{\boldmath $q$}}

\newtheorem{theorem}{Theorem}[section]
\newtheorem{prop}[theorem]{Proposition}
\newtheorem{exam}[theorem]{Example}

\newtheorem{defs}[theorem]{Definition}
\newtheorem{rem}[theorem]{Remark}

\newtheorem{cor}[theorem]{Corollary}
\newtheorem{assump}[theorem]{Assumptions}
\newtheorem{notation}[theorem]{Notation}

\begin{document}

\maketitle

\abstract
\begin{small}
For a   smooth surface in $\mathbb{R}^3$ this article contains local study of certain affine equidistants, that is loci of points at a fixed ratio between
 points of contact of parallel tangent planes (but excluding ratios 0 and 1 where the equidistant contains one or other
point of contact). The situation studied occurs generically in a 1-parameter family, where two parabolic points of the surface have
parallel tangent planes at which the unique  asymptotic directions are also parallel. The singularities
are classified by regarding the equidistants as  critical values of a  2-parameter unfolding of maps from $\mathbb{R}^4$ to $\mathbb{R}^3$.
In particular, the singularities that occur near the so-called `supercaustic chord', joining the two special parabolic points, are classified. For a
given ratio along this chord either one or three special points are identified at which singularities of the equidistant become more special. Many of the
resulting singularities have occurred before in the literature in abstract classifications, so the article also provides a natural geometric
setting for these singularities, relating back to the geometry of the surfaces from which they are derived.
\end{small}

\noindent
MR Classification 57R45, 53A05


\section{Introduction}
 A smooth closed surface in affine 3-space will contain pairs of points at which the affine tangent planes are parallel; indeed
 the tangent plane at a given point may be parallel to that at several other points if the surface is non-convex. Associated with these
 pairs of points, and the chords joining them, there are a number of affinely invariant
 constructions.  The {\em affine equidistants} are the loci of points at a fixed ratio $\lambda : 1-\lambda$ along the chords, and the
 {\em centre symmetry set} is the envelope of the chords, which can be locally empty.  These constructions have been
 examined from the point of view of singularity theory in the last few years by several authors; there are many connexions with
 earlier work such as the `Wigner caustic' of Berry~\cite{berry} which, for a curve in the plane, is the equidistant corresponding
 to a ratio $\lambda=\frac{1}{2}$, that is the midpoints of the parallel tangent chords, and the bifurcations of central symmetry of Janeczko~\cite{janeczko}.
 Notable among recent studies is the work of Domitrz and his co-authors, for example~\cite{domitrz1}.

 A generic surface $M$ in affine 3-space will generically have pairs of points at which the tangent planes are parallel and for which
 both points in the pair are parabolic points of $M$.  For the locus of parabolic points of $M$ is generically a 1-dimensional
 set, a union of smooth curves, and requiring parallel tangent planes imposes  two conditions on a pair of points of this set,
 so that a finite number of solutions can be expected. In this article we investigate one possible local degeneration of this
 generic situation by requiring  also that the unique asymptotic directions coincide
 at such a pair of parabolic points with parallel tangent
 planes. For this to occur the surface $M$ must be contained in a smoothly varying family $M_\ee$ of surfaces. Since
 our  investigation is local we shall in fact consider two surface patches $M_0$ and $N_0$ which vary in a 1-parameter family
 $M_\ee, N_\ee$.  A similar degeneracy was investigated for plane curves in~\cite{GR2017}; we sometimes call it a `supercaustic'
 situation. This term is defined in \S\ref{ss:superc}.

We find the values $\lambda \ne 0,1$ for which the ratio $\lambda : 1-\lambda$
determines an equidistant at which the structure undergoes a qualitative change. There are one or three of these
values, depending on the relative orientation of $M_0$ and $N_0$. One `degenerate' value always exists and
results in a high codimension singularity; we are able to give a partial analysis of this case. When the other
two values exist we call them {\em special values} (Definition~\ref{def:special}), and a complete
analysis is given.

 The article is organized as follows.  In \S\ref{s:general} we introduce the family of surfaces we shall work with (\S\ref{ss:surfaces}),
 and the maps which we shall classify up to $\mathcal A$-equivalence to study the equidistants (\S\S\ref{ss:generating},
 \ref{ss:superc}).  We also show how some of the conditions that arise later can be interpreted geometrically in terms
 of a scaled reflexion map (\S\ref{ss:contact}, Definition~\ref{def:scaledcontactmap}).
 In \S\ref{s:normal} we find normal forms of maps up to $\mathcal A$-equivalence
 that generate the equidistants: they are the sets of critical values of these maps. We examine in that section general
 values of the ratio (Generic Case 1.1)  and the two `special' values (Special Case 1.2), leaving the `degenerate' value (Degenerate Case 2) to \S\ref{s:degen}.

 The main results are contained in Proposition~\ref{prop:def-indef} and the accompanying Figure~\ref{fig:def-indef} for Generic
Case~1.1; Proposition~\ref{prop:case1.2} and the accompanying Figure~\ref{fig:special-clock} for Special Case~1.2,  and  Table~\ref{table1} in~\S\ref{ss:examples} for Degenerate
Case~2.

\section{The general setup}\label{s:general}
\subsection{A generic family of surfaces}\label{ss:surfaces}
Consider the parabolic set $P$ (assumed to be a nonempty smooth curve) of a generic smooth closed surface $M$  in $\RR^3$.  We can expect generically to find
a finite number of  pairs of distinct points on $P$ for which the
tangent planes to $M$ are parallel,  since the two points give us two degrees of freedom and it is two conditions for the tangent planes to be parallel. However it will not be
generically true that the unique asymptotic directions at such a pair of points are parallel. For that we require a 1-parameter family of surfaces and it is this situation which
we study here.

Our considerations are local, and also affinely invariant. For this situation we have two surfaces, $M_\ee$ and $N_\ee$, varying in a 1-parameter family;
using a family of affine transformations of $\RR^3$
(coordinates $(x,y,z)$) we can assume
that the origin lies on $M_\ee$,  that the origin is a parabolic point of $M_\ee$ and that the unique asymptotic direction there is always along the $y$-axis, for all $\ee$ close to 0.
Further we can assume that the point $(0,0,1)$ lies on $N_\ee$ for all small $\ee$ and that for $\ee=0$  this point is parabolic, has horizontal tangent plane parallel
to the $(x,y)$-plane,  and has unique asymptotic direction parallel to the
$y$-axis.   We realise this setup by the surfaces
\begin{eqnarray}
M_\ee:  z = f(x,y,\ee) &=& f_{20}x^2 + f_{300}x^3+f_{210}x^2y + f_{120}xy^2 + f_{030}y^3 + \ldots \nonumber \\
&+& \ee\left( f_{301}x^3+f_{211}x^2y+\ldots\right) + \ee^2\left(f_{302}x^3 + \ldots \right) + \ldots, \label{eq:M}\\
N_\ee: z = 1 + g(x,y,\ee) &=& 1 + g_{20}x^2 +   g_{300}x^3+g_{210}x^2y + g_{120}xy^2 + g_{030}y^3 + \ldots \nonumber \\
&+& \ee\left(g_{101}x+g_{011}y + g_{201}x^2 + g_{111}xy+g_{021}y^2 + \ldots\right)  \nonumber \\
&+& \ee^2\left( g_{102}x + g_{012}y + \ldots\right) + \ldots. \label{eq:N}
\end{eqnarray}
For terms other than $f_{20}, g_{20}$, subscripts $ijk$ indicate that the corresponding monomial is $\ee^kx^iy^j$.

We make the following assumptions about these expansions.
\begin{assump}
{\rm (i) $f_{20}\ne 0, g_{20}\ne 0$, that is neither $M_0$ nor $N_0$ is umbilic at its basepoint $(0,0,0)$ or $(0,0,1)$.

(ii) $f_{030} \ne 0, g_{030}\ne 0$, that is the parabolic curves of $M_0$ at the origin and $N_0$ at $(0,0,1)$ are smooth and not tangent
to the asymptotic directions there (i.e.\ these points are not cusps of Gauss).  We shall take $f_{030}>0$ without loss of generality,
and we sometimes write $f_{030}=f_3^2, \ g_{030}=\pm g_3^2$ when a definite sign is needed, to avoid square roots appearing in the formulas.
}
\label{assump}
\end{assump}

\subsection{Family of maps for the equidistants}\label{ss:generating}
The $\lambda$-equidistant for a fixed $\ee$ is the locus of points in $\RR^3$ of the form $(1-\lambda)\bp \,+\lambda\bq$ where $\bp\in M_\ee, \bq\in N_\ee$ and
the tangent planes to $M_\ee$ at $\bp$ and $N_\ee$ at $\bq$ are parallel.

\smallskip

{\em We always assume $\lambda \ne 0, \lambda \ne 1$ in what follows.}

\smallskip

We use $s=(s_1, s_2)$ as parameters on $M_\ee$ and similarly $t=(t_1,t_2)$ for $N_\ee$; we have a 2-parameter family of maps $\RR^4\to\RR^3$:
\begin{equation}
 \RR^4 \times \RR^2 \to \RR^3, \ (s,t, \ee,\lambda) \mapsto (1-\lambda)(s_1,s_2,f(s_1, s_2, \ee)) + \lambda(t_1, t_2, 1+g(t_1, t_2, \ee)).
\label{eq:H1}
\end{equation}
Then it is straightforward to check that, for fixed $\ee$ and $\gl$, the set of critical values of this map is the $\lambda$-equidistant of $M_\ee$ and $N_\ee$.  We are therefore
interested in this family of maps up to $\cA$-equivalence.  We make the
change of variables
\[ (1-\lambda)s_1+\gl t_1 = u_1, \ (1-\lambda)s_2 + \gl t_2 = u_2, \mbox{ and write } \lambda=\lambda_0+\alpha, \]
replacing $t_1$ and $t_2$, to rewrite (\ref{eq:H1}) as a map of the form (for any $\lambda_0\ne 0,1$)
\begin{equation}
H: \RR^4\times\RR^2 \to \RR^3, \ H(s_1,s_2,u_1,u_2,\ee,\alpha)= (u_1, u_2, h(s_1,s_2,u_1,u_2,\ee,\gl_0+\alpha)).
\label{eq:H}
\end{equation}
regarded as a 2-parameter unfolding of the map $H_0(s_1,s_2,u_1,u_2,0,\gl_0)$. Therefore we have the following.
\begin{prop}
The $\lambda$-equidistant for fixed $\ee$ is the set of points $(u_1, u_2, h)\in \RR^3$ for which $\partial h/\partial s_1 = \partial h/\partial s_2 = 0.$
For fixed $\gl$ the union of all the equidistants, spread out in $\RR^4$,
the  planar sections of which are the $\ee=$ constant equidistants, is the set of points $ (u_1, u_2, h, \ee) \in\RR^4$  for
which the same conditions
$\partial h/\partial s_1 = \partial h/\partial s_2 = 0$ hold.
\end{prop}

\subsection{Maps and supercaustics}\label{ss:superc}
Let $\phi:\RR^4\to \RR^2$ be given, for fixed $\lambda$ and $\ee$,  by $\phi(s_1,s_2,u_1,u_2)=(h_{s_1}, h_{s_2})$, subscripts denoting partial derivatives as usual.
Then the corresponding equidistant, given by $\phi^{-1}(0,0)$, is singular when there is a kernel vector of $d\phi$ with image under $dH$ equal to {\bf 0}, these being evaluated
at a point of $\phi^{-1}(0,0)$. This requires that
\[ \mbox{rank } J < 4 \mbox{ where } J = \left(\begin{array}{cccc}h_{s_1s_1}&h_{s_1s_2}&h_{s_1u_1}&h_{s_1u_2}\\
h_{s_2s_1}&h_{s_2s_2}&h_{s_2u_1}&h_{s_2u_2}\\
0&0&h_{u_1}&h_{u_2}\\
0&0&1&0\\
0&0&0&1 \end{array}\right),\]
that is $h_{s_1s_1}h_{s_2s_2}=h_{s_1s_2}^2.$  The singular points of the equidistant for fixed $\lambda$ and $\ee$ are therefore
\begin{equation}
 \{ (u_1, u_2, h(s_1,s_2,u_1,u_2))  \  : \  h_{s_1}=h_{s_2}=h_{s_1s_1}h_{s_2s_2}-h_{s_1s_2}^2=0\}.
\label{eq:singpts}
\end{equation}
We note here that, for fixed $\ee$, the `centre symmetry set' of the pair of surfaces $M,N$ \cite{GZ}, which is the locus of singular points of the equidistants for varying $\lambda$,
is given by the same formula (\ref{eq:singpts}) where $h$ is now a function of $s_1, s_2, u_1, u_2, \lambda$ but with $\ee$ still fixed.

It is possible that some singular points of the equidistant arise from singularities of the critical set itself in $\RR^4$. In our case this requires, for fixed $\lambda$ and $\ee$,  that
the top two rows of the above matrix $J$ are dependent.  Indeed, evaluating these rows at $(s_1,s_2,u_1,u_2,\lambda,\ee)=(0,0,0,0,\lambda,0)$ the second row is entirely zero.
This means that, for all $\lambda$, but $\ee=0$, the critical set itself is singular at the origin of $\RR^4$.
\begin{defs}
{\rm
In the above situation, the $\lambda$-axis is called a {\em supercaustic}; see \cite{GR2017}. The whole of this axis
maps to singular points of the equidistants.
}
\end{defs}

\begin{rem}
{\rm
This depends crucially on the special nature of our surfaces, with not only parallel tangent planes at parabolic points of $M_0$ and $N_0$ but also the asymptotic directions at those
points being parallel.  If instead we assume that the asymptotic directions are distinct (without loss of generality
we can take them along the $x$ and $y$ axes) then the top two rows of $J$ become independent for
$s_1=s_2=u_1=u_2=\ee=0$ and arbitrary $\lambda$. In fact, writing $g_{020}$ for the coefficient of $y^2$ in the parametrization of $N_0$  and putting $g_{20}=0$ these rows
become
\[ \left(\begin{array}{cccc} 2(1-\lambda)f_{20} & 0 & 0 & 0 \\
0 & \frac{2g_{020}(1-\lambda)^2}{\lambda} &0 & -\frac{2g_{020}(1-\lambda)}{\lambda} \end{array}\right). \]
In this case the `supercaustic' is empty.
}
\end{rem}

\subsection{Scaled reflexion map and contact} \label{ss:contact}

Consider the affine map $\cS: \RR^3\to\RR^3$ given by $\cS(x,y,z) = (\mu x, \mu y, \mu(z-1))$ where $\mu=\frac{\lambda}{\lambda-1}\ne 0$.  This leaves the point $(0,0,\lambda)$ fixed
and maps  $(0,0,1)$ to the origin.  We can measure the contact between $\cS(N_0)$ and $M_0$ by
composing  the parametrization of $\cS(N_0)$
given by $\left(\mu x, \mu y, \mu g(x, y, 0)\right) $ with the equation of $M_0$, say $Z-f(X, Y, 0)=0$.
\begin{defs}{\rm
The {\em scaled contact map} is the contact map germ
}
\[ K:\RR^2, (0,0) \to \RR, 0, \ K(x,y)=\mu g(x,y,0) - f(\mu x, \mu y, 0), \ \mu=\frac{\lambda}{\lambda-1} \ \ \ {\rm as \ above}. \]
\label{def:scaledcontactmap}
\end{defs}
\vspace*{-0.7cm}
We  shall find this contact map useful in interpreting the conditions which arise from $\ee$-families of equidistants as $\ee$ passes through 0.

The 2-jet of $K$ is $K_2(x,y)=\mu(g_{20}-\mu f_{20})x^2$ so that in our situation $K$ is always non-Morse; it has corank 1 and is of
type $A_k$ at $(0,0)$ for some $k$, provided $f_{20}\lambda + g_{20}(1-\lambda)\ne 0$
(when this fails we call this the `Degenerate Case~2'; see \S\ref{s:degen}).
The coefficient of $y^3$ in $K$ is $\mu(g_{030}-\mu^2f_{030})$ so that $K$ is then of type exactly $A_3$ provided $f_{030}\lambda^2 - g_{030}(1-\lambda)^2 \ne 0$. If $f_{030}, g_{030}$  are nonzero and have
opposite signs then of course this coefficient can never be zero.

\begin{defs}{\rm
Assume as above that $f_{20}\lambda + g_{20}(1-\lambda)\ne 0$.
When $f_{030}, g_{030}$ have the same sign (without loss of generality, positive), and the above coefficient  $f_{030}\lambda^2 - g_{030}(1-\lambda)^2$ of $y^3$ is zero, then we refer to the two resulting
values of $\lambda$ as {\em special values}.  Writing $f_{030}=f_3^2, g_{030}=g_3^2$ where we may take $f_3>0, g_3>0$,
these special values of $\lambda$ are $\displaystyle{\frac{g_3}{g_3\pm  f_3}}$.  (We shall usually assume $f_3 \ne g_3$ to avoid one of the special values `going to infinity'.)
These special values of $\lambda$ give rise to what we shall call Special Case~1.2. This is examined in detail in \S\ref{ss:special}.
\label{def:special}
}
\end{defs}

When $\lambda$ has a special value, say $\displaystyle{\frac{g_3}{g_3 + f_3}}$, the condition for $K$ to have exactly type $A_3$
at $(0,0)$ works out to be
\begin{equation}
 (4g_{040}g_{20}-g_{120}^2)f_3^4 + 4g_{040}f_{20}f_3^3g_3+2f_{120}g_{120}f_3^2g_3^2+4f_{040}g_{20}f_3g_3^3 + (4f_{040}f_{20}-f_{120}^2)g_3^4 \ne 0.
\label{eq:A3}
\end{equation}
This condition will be satisfied by a generic pair of surfaces $M_0, N_0$.
With the other special value the signs in front of the coefficients of $f_3^3g_3$ and $f_3g_3^3$ both change to minus.

When the quadratic terms of the contact map $K$ vanish identically, that is when $f_{20}\lambda + g_{20}(1-\lambda)= 0$, the cubic terms will in general be nondegenerate and $K$ will
generically have type $D_4^\pm$, that is $\mathcal{R}$-equivalent to $x^3\pm xy^2$. The polynomial in the coefficients of $f$ and $g$ which distinguishes the two cases is rather complicated
but, remarkably, it has a different interpretation which we give in \S\ref{s:degen} in the context of self-intersections of the equidistant.  See Remark~\S\ref{rem:D4-contact}.

\section{The equidistants: normal forms}\label{s:normal}
For a general study of the equidistants we need to expand the function $h$ in (\ref{eq:H}) using the parametrizations (\ref{eq:M}) and (\ref{eq:N}).
We begin with $\ee=0$ and write, for a fixed $\gl$, $H_{0\gl}(s,u)=(u, h_{0\gl}(s,u))=H(s,u,0,\gl)$.
The coefficient of $s_1^is_2^ju_1^ku_2^\ell$ in $h_{0\gl}$ will be written $c_{ijk\ell}$.
We find:
\[ \mbox{The 2-jet of } h_{0\gl} \mbox{ at  } s=u=0 \mbox{ is } (1-\gl)(\gl f_{20} + (1-\gl)g_{20})s_1^2 - 2g_{20}\textstyle{\frac{1-\gl}{\gl}}s_1u_1.\]
Note that the coefficient of $s_1u_1$ is nonzero.

The main subdivision is between those $\gl$ for what  $\gl f_{20} + (1-\gl)g_{20}$ is nonzero (Generic Case~1) or zero (Degenerate Case 2).  We cover the Generic Case
here and the Degenerate Case in \S\ref{s:degen} below.

\medskip\noindent
{\bf Case 1} \ $\gl f_{20} + (1-\gl)g_{20} \ne 0$.  From \S\ref{ss:contact} this is also the condition for the contact function $K$ to have type
$A_k$ for some $k$.

We can now redefine the variable $s_1$ (`completing the square') to eliminate all terms containing $s_1$ besides $s_1^2$ in $h_{0\gl}$. The coefficient of $s_2^3$ then becomes
\[ c_{0300}=\textstyle{\frac{1-\gl}{\gl^2}}\displaystyle (f_{030}\gl^2-g_{030}(1-\gl)^2).\]

\subsection{The general values of $\lambda$}
{\bf Generic Case 1.1} $c_{0300}\ne 0$, that is, $Q \ne 0$ where
 \begin{equation}
 Q = f_{030}\gl^2-g_{030}(1-\gl)^2.
 \label{eq:Q}
 \end{equation}
 From \S\ref{ss:contact} this is also the condition for the contact function $K$ to have type $A_2$ and that $\gl$ is not a {\em special value}.

Consider the 3-jet of $H_{0\gl}$. There are six degree 3 monomials which do not involve $s_1$ and which do involve $s_2$ (any monomial in $u_1, u_2$ alone can
be eliminated by a `left-change' of coordinates). We still have the freedom to change coordinates in $s_2$ (involving $s_2, u_1, u_2$) and in $u_1, u_2$ (involving
$u_1, u_2$ only). Using only the first of these the terms in $s_2^2u_1$ and $s_2^2u_2$ can be eliminated, leaving
\begin{equation}
\left(u_1, u_2, (1-\gl)(\gl f_{20} + (1-\gl)g_{20})s_1^2 +c_{0300}s_2^3 + s_2\left(c_{0120}u_1^2+c_{0111}u_1u_2+c_{0102}u_2^2\right)\right).
\label{eq:3jet}
\end{equation}
(The coefficients $c_{ijk\ell}$  need to be updated to take account of the substitutions.)
The quadratic form in $u_1$ and $u_2$ can be diagonalised, eliminating the term in $s_2u_1u_2$ so that, scaling $s_1$, the
last coordinate in $\RR^3$ and $s_2$, we have 3-jet, say
\[ (u_1, u_2, s_1^2+s_2^3+as_2u_1^2+bs_2u_2^2).\]

Suppose that the quadratic form in parentheses in (\ref{eq:3jet}) is not a perfect square, that is $c_{0111}^2-4c_{0120}c_{0102}\ne 0$.
Then $a$ and $b$ above are nonzero.
The condition for this  is $R\ne 0$ where
\begin{equation}
R = f_{20}^2f_{030}\left(g_{120}^2-3g_{210}g_{030}\right) - g_{20}^2g_{030}\left(f_{120}^2-3f_{210}f_{030}\right).
\label{eq:R}
\end{equation}
Since this condition does not involve $\gl$ it will  be satisfied by a generic pair of surfaces $M_0, N_0$. Note that the condition separates into a quantity for $M_0$ unequal to the same quantity for $N_0$.

\begin{prop}
The condition $R\ne 0$ can also be interpreted as saying that the images under the Gauss map of the parabolic curves on $M_0$ and $N_0$
have ordinary tangency (that is, 2-point contact) in the Gauss sphere.  These images are smooth by Assumptions~\ref{assump}.
\label{rem:contact}
\end{prop}
{\bf Proof} \ The parabolic curves on the two surfaces are given by $f_{xx}f_{yy}-f_{xy}^2=0$ and $g_{xx}g_{yy}-g_{xy}^2=0$
for $M_0$ and $N_0$ respectively. The surface $M_0$ has a parabolic point at the origin and $N_0$ has a parabolic point
at $(0,0,1)$ and since they have parallel asymptotic directions at these points the images of the
respective parabolic curves under the Gauss map are
tangent.  We shall use the modified Gauss maps, that is $(x,y)\mapsto (X,Y)=(f_x, f_y)$ and similarly for $g$.  By a direct
calculation, for $M_0$ the image of the
parabolic curve, parametrized by $x$, under the modified Gauss map has an equation, up to terms in $X^2$, of the form
\[ Y=\frac{3f_{030}f_{210}-f_{120}^2}{12f_{20}^2f_{030}}X^2 \]
with a similar result for $N_0$.   The coefficients of $X^2$ are unequal, that is the images have ordinary tangency,
if and only if the condition $R$ above is nonzero. \hfill$\Box$

\medskip

\begin{figure}[!ht]
\begin{center}
\hspace*{-1.5cm}
\scalebox{1.8}{\includegraphics[width=1.3in]{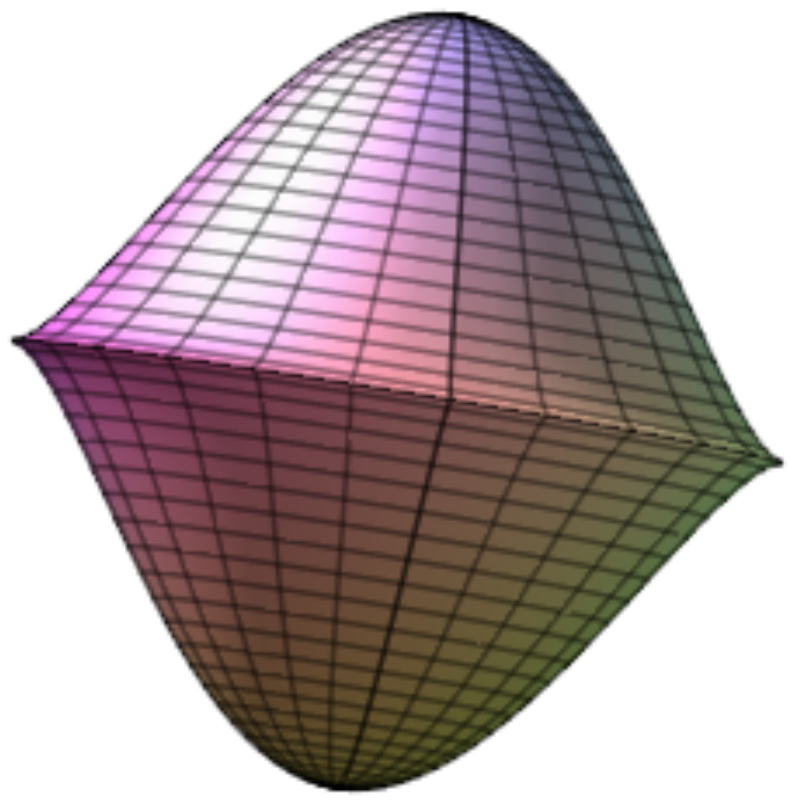}}
\hspace*{-2.5cm}
\scalebox{1.8}{\includegraphics[width=1.3in]{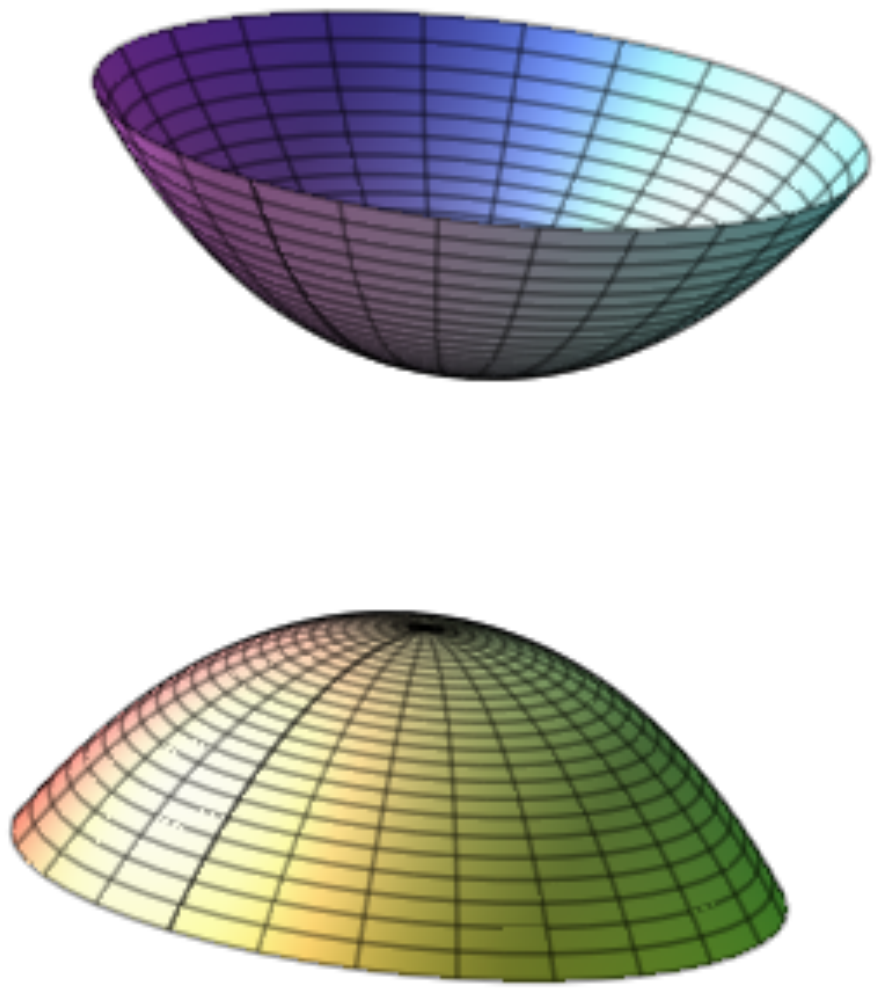}}
\hspace*{-2.5cm}
\scalebox{1.8}{\includegraphics[width=1.3in]{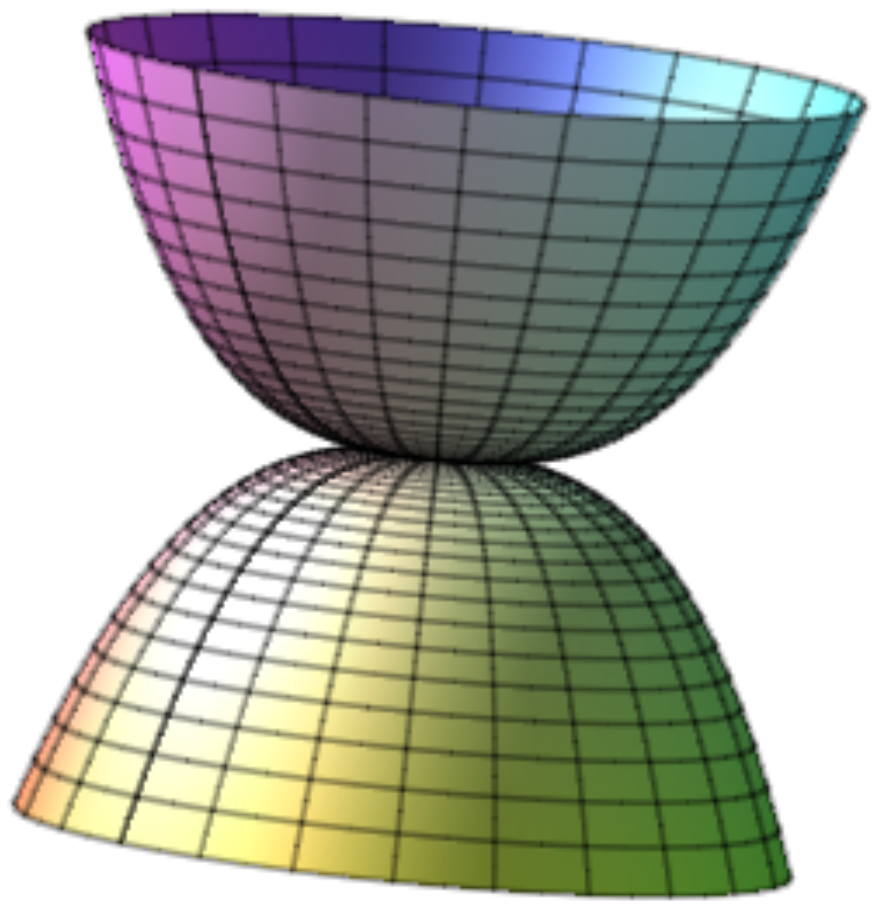}}
\hspace*{-2.7cm}
\scalebox{1.8}{\includegraphics[width=1.3in]{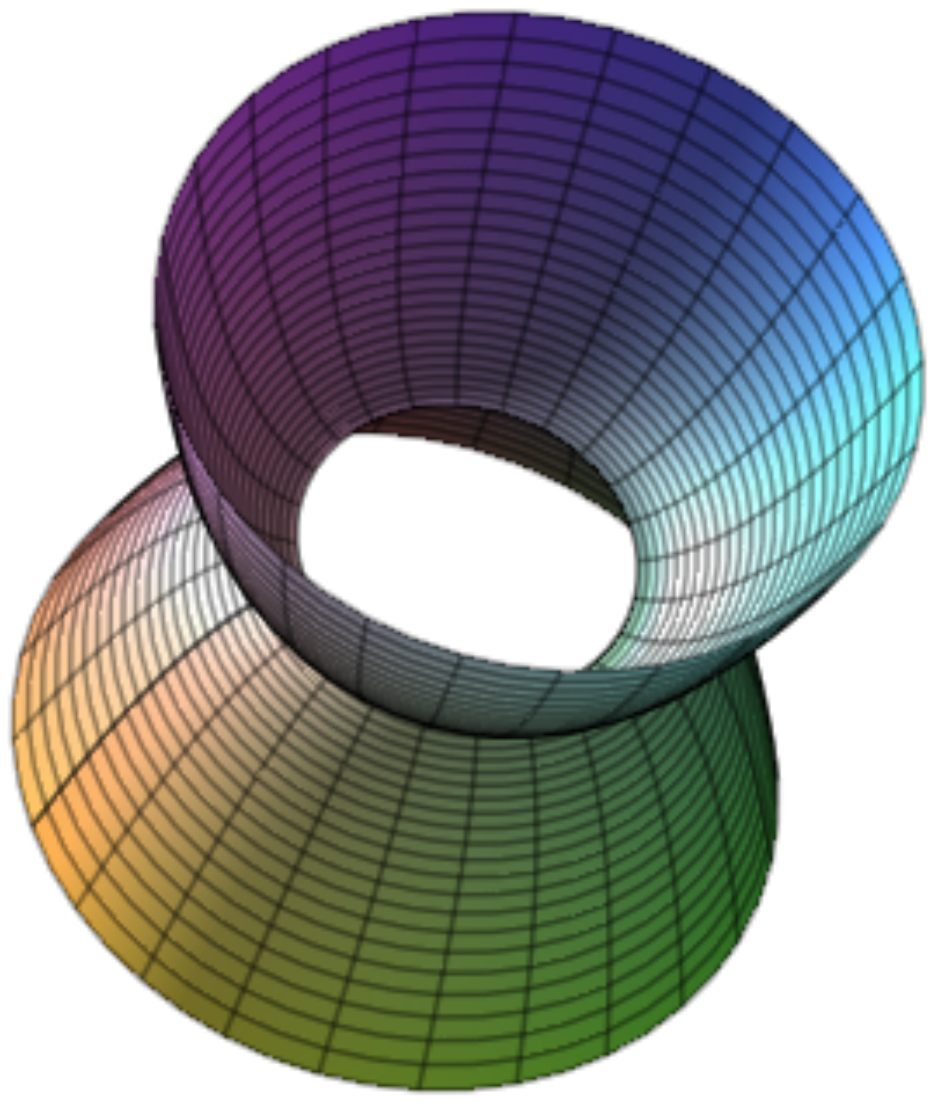}}
\end{center}
\vspace*{-3.5cm}
{\small
\hspace*{1cm}  positive def., $\ee<0$ \hspace{0.2cm} negative def., $\ee<0$ \hspace{0.25cm} negative def., $\ee=0$ \hspace{0.25cm}negative  def., $\ee>0$ }
\vspace*{-0.5cm}
\begin{center}
\scalebox{1.8}{\includegraphics[width=1.3in]{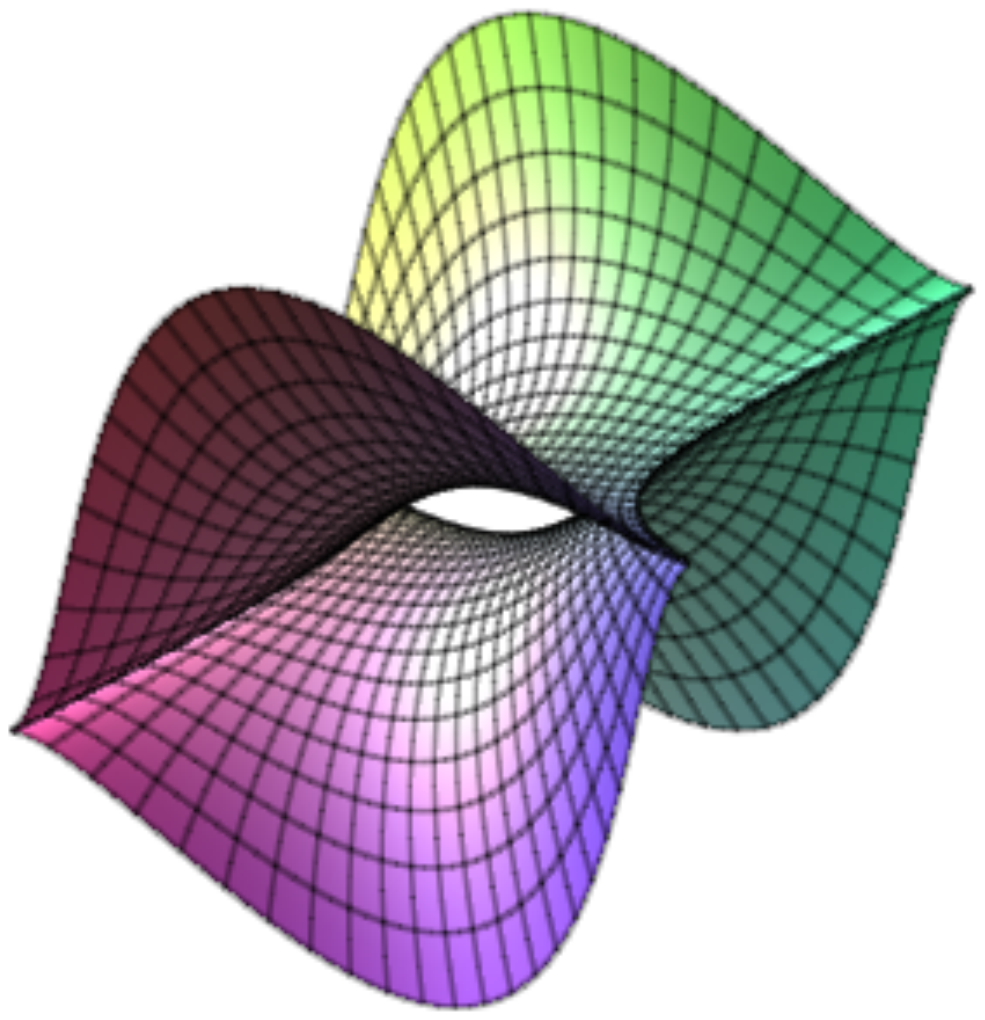}}
\hspace*{-3cm}
\scalebox{1.8}{\includegraphics[width=1.3in]{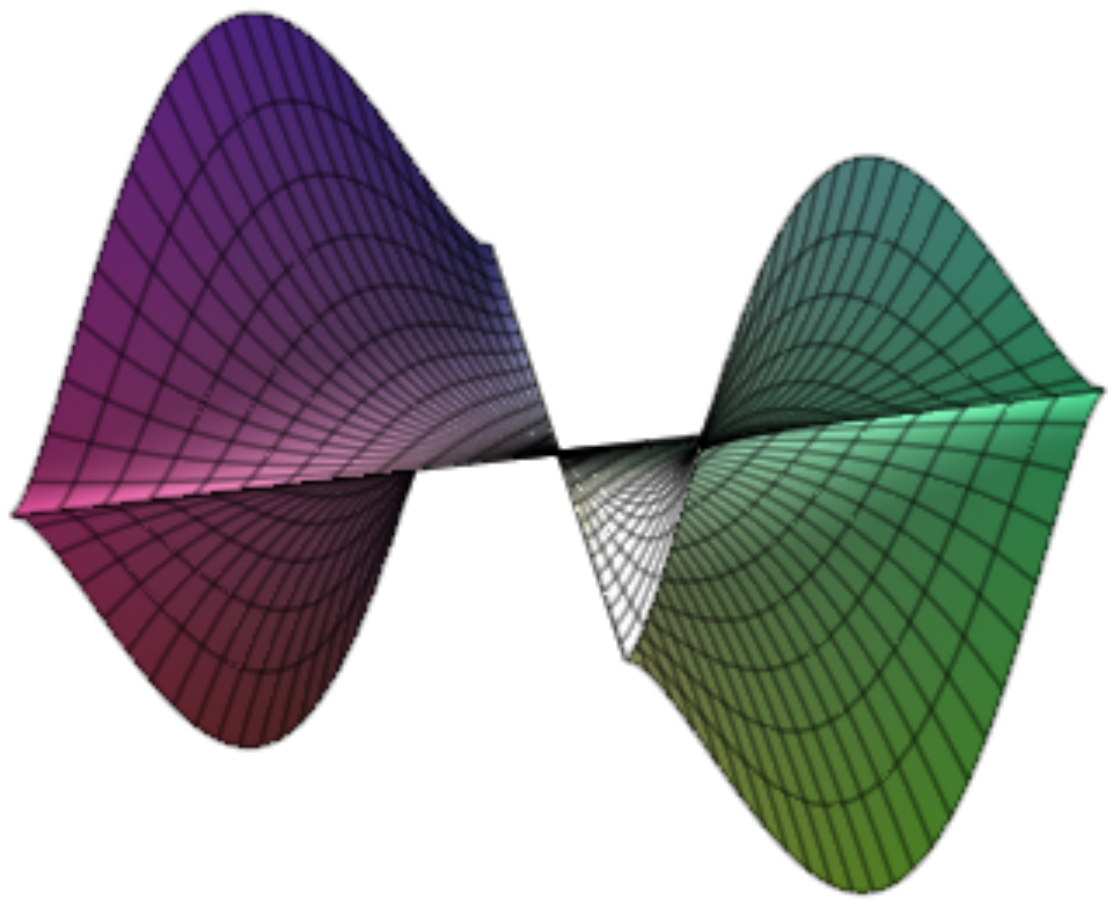}}
\hspace*{-3cm}
\scalebox{1.8}{\includegraphics[width=1.3in]{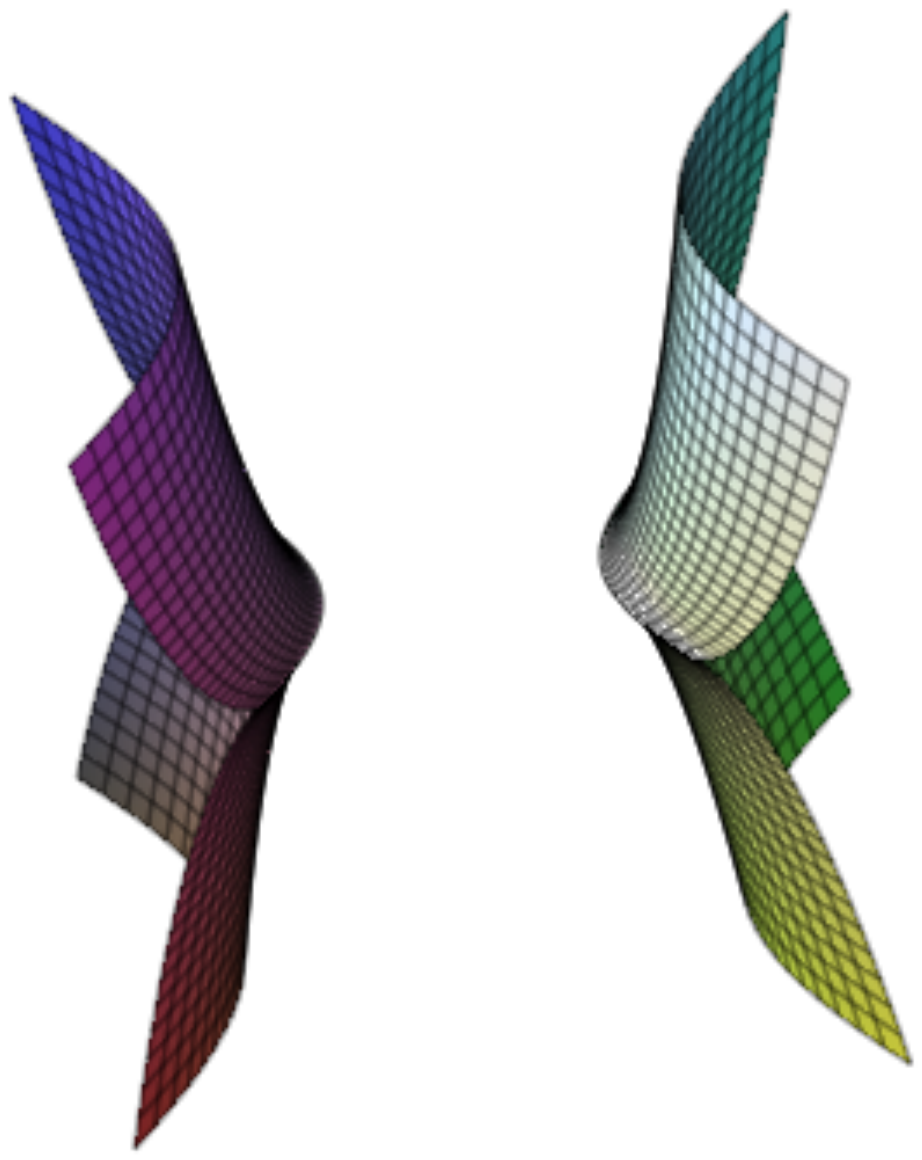}}
\end{center}
\vspace*{-3.5cm}
{\small \hspace{4cm} Indefinite, $\ee<0$ \hspace{0.4cm} Indefinite, $\ee=0$ \hspace{0.4cm} Indefinite, $\ee>0$.}
\caption{\small The various subcases of Proposition \ref{prop:def-indef}:  Positive definite (for $\ee>0$ the equidistant is empty
and for $\ee<0$ has a compact cuspidal edge); 1.1.2 Negative definite,
where for $\ee>0$ there is a compact cuspidal edge; 1.1.3 Indefinite, which has two cuspidal edges
for $\ee\ne 0$ that form a crossing when $\ee=0$. }
\label{fig:def-indef}
\end{figure}

Further scaling allows this case to be reduced to
\begin{equation}
H_{0\gl}(s,u)= (u_1, u_2, s_1^2+s_2^3\pm s_2u_1^2 \pm s_2u_2^2),
\label{eq:1.1}
\end{equation}
where the $\pm$ signs are independent, but by interchanging $u_1$ and $u_2$ we reduce to three cases, as follows.
\begin{prop}
The normal form {\rm (\ref{eq:1.1})} is as follows, using the notation of {\rm (\ref{eq:Q})} and {\rm (\ref{eq:R})}.
See Figure~\ref{fig:def-indef}.

\noindent
{\rm Subcase 1.1.1} (positive definite): $H_{0\gl}(s,u)= (u_1, u_2, s_1^2+s_2^3 +s_2u_1^2 + s_2u_2^2)$.  \\
The condition for this is $f_{030}g_{030} < 0$ and $QR>0$. Bearing in mind the assumptions \ref{assump}
the latter condition is equivalent to $R>0$. This subcase will also be referred to as $A_2^{++}$.\\
{\rm Subcase 1.1.2} (negative definite): $H_{0\gl}(s,u)= (u_1, u_2, s_1^2+s_2^3 -s_2u_1^2 - s_2u_2^2)$.  \\
The condition for this is $f_{030}g_{030} > 0$ and $QR>0$. This subcase will also be referred to as $A_2^{- -}$\\
{\rm Subcase 1.1.3} (indefinite): $H_{0\gl}(s,u)= (u_1, u_2, s_1^2+s_2^3 +s_2u_1^2 - s_2u_2^2)$. \\
The condition for this  is $QR< 0$. In the case when $f_{030}g_{030}<0$ the condition becomes $R<0$.
This subcase will also be referred to as $A_2^{+-}$, \hfill$\Box$
\label{prop:def-indef}
\end{prop}

\smallskip\noindent
The values of $f_{030}, g_{030}$ and $R$ are fixed by the two surfaces $M_0$ and $N_0$.
However, assuming $f_{030}g_{030}>0$, special values of $\lambda$ exist at which $Q$ as in (\ref{eq:Q}) is zero.  Then, as $\lambda$ passes through
such a special value, the normal form changes
between negative definite and indefinite, so that the family of equidistants, for $\ee$ passing through 0, changes accordingly.

\smallskip\noindent
Using standard techniques it can be checked that (\ref{eq:1.1}) is 3-$\cA$-determined, and that an $\cA_e$-versal unfolding is
given by adding a multiple of $(0,0,s_2)$ to the above normal form:
\begin{equation}
H_{\ee\gl}(s,u)= (u_1, u_2, s_1^2+s_2^3\pm s_2u_1^2 \pm s_2u_2^2+\ee s_2).
\label{eq:normalform-unf}
\end{equation}
In terms of the original surfaces the coefficient of $\ee s_2$ is $-g_{011}(1-\lambda)$, and therefore we require
$g_{011}\ne 0$ for a versal unfolding by the parameter $\ee$.

\begin{rem}
{\rm
It is interesting to  relate the above classification to that of the regions on $M$ and $N$ which contribute to the pairs of parallel tangent planes (compare
Prop.2.4 and Figure~3 of \cite{GR2015}). A schematic diagram of the common regions for $M$ and $N$  on the Gauss sphere is  given in Figure~\ref{fig:parab-curves}
below.  The relationship between these and the classification of Proposition~\ref{prop:def-indef} is as follows.

\noindent
Subcase 1.1.1 (positive definite, $f_{030}g_{030} < 0$ and $R>0$): This is (d).\\
Subcase 1.1.2 (negative definite, $f_{030}g_{030} > 0$ and $QR>0$): This is (ac).\\
Subcase 1.1.3 (indefinite): This can arise in two ways, as either (ac) or (b)\\
\hspace*{1cm} (ac) when $f_{030}g_{030}>0$ and $QR<0$,\\
\hspace*{1cm} (b) when $f_{030}g_{030}<0$ and $R < 0$.
}
\end{rem}
\begin{figure}[!ht]
\begin{center}
\scalebox{1.4}{\includegraphics[width=3in]{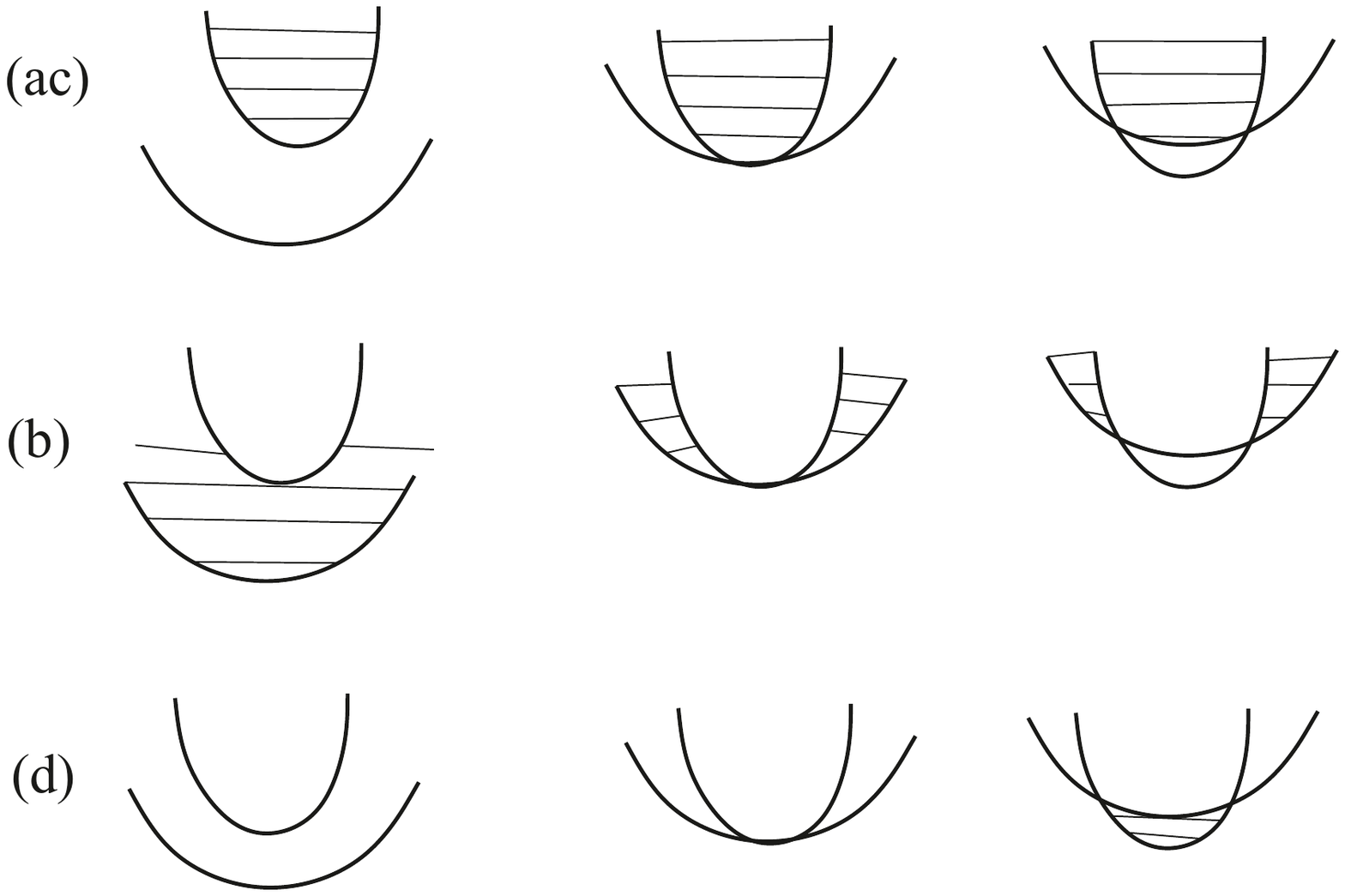}}
\end{center}
\vspace*{-9.5cm}
\caption{\small Schematic diagrams of the images of the Gauss map for the surfaces $M_\ee$ and $N_\ee$.  The curves represent the parabolic curves of these surfaces,
along which the Gauss map has a fold,
and the hatched regions represent the regions where the images of the Gauss maps of
$M_\ee$ and $N_\ee$ intersect, that is the regions of the Gauss sphere representing parallel normals (or parallel tangent planes). Left to right of each row shows varying $\ee$, with the
middle diagram $\ee=0$, and the three possible cases are labelled (ac), (b), (d) as described in the text, to accord with Figure~3 in \cite{GR2015}. Note that
the two curves for $\ee=0$ have ordinary tangency---see Remark~\ref{rem:contact}. }
\label{fig:parab-curves}
\end{figure}

\begin{figure}[!ht]
\begin{center}
\scalebox{1.4}{\includegraphics[width=4in]{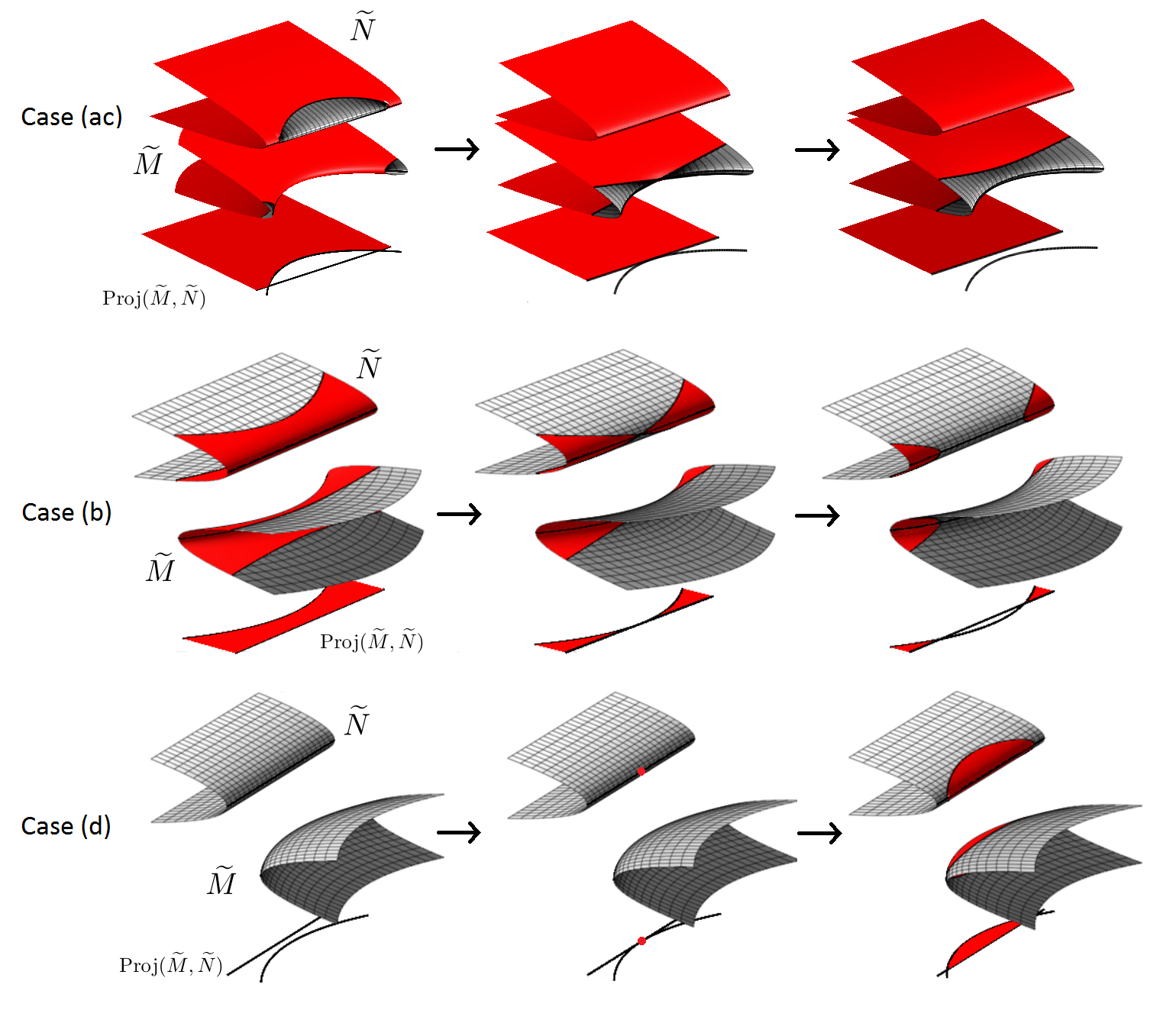}}
\end{center}
\vspace*{-0.5cm}
\caption{\small In this  diagram, the Gauss map of the surfaces $M_\ee$ and $N_\ee$ is represented by vertical projection and the surfaces
in this schematic representation are labelled $\widetilde{M}, \widetilde{N}$.  The rows and columns are arranged as in Figure~\ref{fig:parab-curves}. See the above text
for further explanation. }
\label{fig:Gauss-Map-Projections-3cases}
\end{figure}

Let us call a pair of points, one from $M_\ee$ and the other from $N_\ee$, at which the tangent planes are parallel, `mates'.
Consider for example the top left diagram of Figure~\ref{fig:parab-curves} and assume that the upper curve is the image of the parabolic curve of $N_\ee$ in the Gauss sphere.
Each point above this curve is the image of two points of $N_\ee$ and two points of $M_\ee$ giving altogether four mates. Each point on the upper curve
is the image of two points of $M_\ee$ and a single parabolic point of $N_\ee$ which is a mate for both of them.  On the surface $M_\ee$ itself
there will be a region close to the base-point $(0,0,0)$ consisting of those points of  $M_\ee$ with at least one mate, and usually two mates, on $N_\ee$---a region
`doubly covered by mates on $N_\ee$'. This region
will have a local boundary corresponding in the way just described to the parabolic curve on $N_\ee$. Turning to the upper right diagram of Figure~\ref{fig:parab-curves}
the hatched region representing mates  now contains a segment of the parabolic curve of $M_\ee$. On the surface $N_\ee$ this will result in a closed loop
on the boundary of the region of points having mates on $M_\ee$. The situation on the surfaces themselves is illustrated schematically in Figure~\ref{fig:Gauss-Map-Projections-3cases}.

\subsection{The `special values' of $\lambda$}\label{ss:special}
{\bf Special Case 1.2} $c_{0300} = 0$, that is $\gl$ has one of the two special values as in \S\ref{ss:contact}.  Note that this requires $f_{030}$ and $g_{030}$
to have the same sign, which we take as positive, and write $f_{030}=f_3^2, g_{030}=g_3^2$ where $f_3>0, g_3>0$.

This case will be examined by choosing one of the special values for $\gl$ given by $c_{0300}=0$, namely $\gl = \displaystyle{\frac{g_3}{g_3 +  f_3}}$.
We can eliminate the terms in $s_2u_2^2$ and $s_2u_1u_2$ by a substitution of the form $s_2=s'_2 +au_1+bu_2$, assuming only
the condition $\gl f_{20} + (1-\gl)g_{20} \ne 0$ of Generic Case 1. The coefficient of $s_2^2u_2$ then
becomes $3f_2^2\ne 0$ and the remaining degree 3 terms in $h_{0\gl}$, namely $s_2^2u_1, s_2^2u_2$ and $s_2u_1^2$  can therefore be reduced
to the last two by redefining $u_2$, at the same time making the coefficient of $s_2^2u_2$ equal to 1. The 3-jet of $H_{0\gl}$ is now of the form (scaling $s_1$)
\[ (u_1, u_2, \pm s_1^2 +s_2^2u_2+c_{0120}s_2u_1^2),\]
where the updated $c_{0120}$ is nonzero if and only if  $R\ne 0$ as in (\ref{eq:R}), and for generic $M_0, N_0$ this will be satisfied.

Passing to the 4-jet of $H_{0\gl}$, we can first remove all monomials divisible by $s_1$ besides $\pm s_1^2$ by completing the square, and then
eliminate all degree 4 monomials besides $s_2^4$ and $s_2^3u_1$, without
adding any new monomials of degree 3.  This can be done,
for example, by substitutions of the form $s_2 = s'_2 + \mbox{ quadratic terms in } s'_2, u'_1, u'_2$, $u_1 = u'_1 + \mbox{ quadratic terms in } u'_1, u'_2,$
and similarly for $u_2$. A left change of coordinates will then restore the first two components of $H_{0\lambda}$ to $(u_1, u_2)$.

The 4-jet is now reduced to
\[ \left(u_1, u_2, \pm s_1^2+s_2^2u_2+c_{0120}s_2u_1^2+c_{0400}s_2^4+c_{0310}s_2^3u_1\right).\]
This is 4-$\cA$-determined provided all the coefficients are nonzero.
The coefficient $c_{0400}$ is nonzero if and only if the `exactly $A_3$ contact condition' (\ref{eq:A3}) holds.
Unfortunately we do not know a geometrical criterion for the coefficient of $s_2^2u_2$ to be nonzero;
it involves only the coefficients in the functions $f, g$ which define the surfaces $M_0, N_0$.

\begin{figure}
\begin{center}
\scalebox{1.4}{\includegraphics[width=5in]{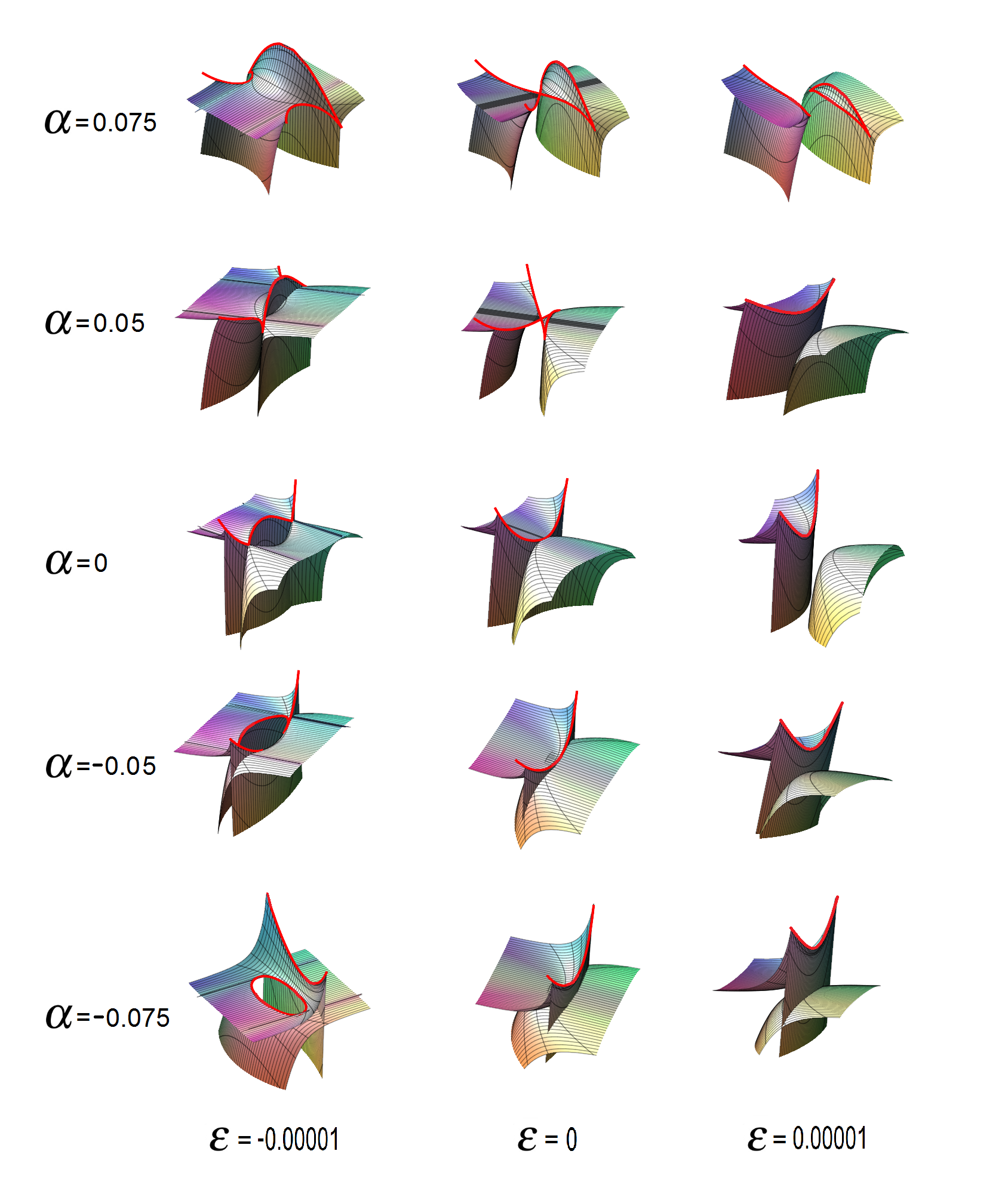}}
\end{center}
\caption{\small Special Case 1.2. A typical `clock diagram' of equidistants close to a special value of $\lambda_0=g_3/(g_3\pm f_3)$. The vertical axis
represents $\lambda=\lambda_0+\alpha$ and the horizontal axis the parameter $\ee$ in the family of surfaces.}
\label{fig:special-clock}
\end{figure}

Scaling reduces all but the coefficient of $s_1^2$ to 1 and we summarize this discussion as follows.
\begin{prop}
For Special Case 1.2, that is $f_{030}=f_3^2, g_{030}=g_3^2$, a special value of $\lambda=g_3/(g_3\pm f_3)$ (Definition~\ref{def:special} or $Q=0$ as in (\ref{eq:Q}))
 but $\gl f_{20} + (1-\gl)g_{20} \ne 0$,
the function $H_{0\lambda}$ reduces under $\mathcal A$-equivalence to the normal form
\begin{equation}
H_{0\lambda}(s_1,s_2,u_1,u_2)=\left(u_1, u_2, \pm s_1^2+s_2^2u_2+s_2u_1^2+ s_2^4+s_2^3u_1 + (ps_2+qs_2^3) \right),
\label{eq:1.2}
\end{equation}
provided the geometrical conditions $R\ne 0$ (\ref{eq:R}), and `exactly $A_3$-contact' (\ref{eq:A3}) hold, together with a third
condition on $M_0, N_0$ which will be generically satisfied.  The terms $ps_2+qs_2^3$ in brackets represent an $\cA_e$-versal unfolding
provided the geometrical condition $g_{011}\ne 0$ in (\ref{eq:M}) holds.  See Figure~\ref{fig:special-clock}
for a `clock diagram' of the equidistants in the $(p,q)$-plane. \hfill$\Box$
\label{prop:case1.2}
\end{prop}
A similar normal form, without the fourth variable $s_1$, but with an additional ambiguity of sign, occurs as $4_2^2$ in \cite{Marar-Tari}; see also  \cite{Gory}.
 The sign in front of $s_1^2$ will not affect our
results since the critical set of $H_{0\lambda}$ has $s_1=0$.
The versal unfolding condition means that as $\ee$ changes through 0 the normal
to $N$ tilts in a direction with a nonzero component along the $y$-axis, which is the asymptotic direction at $\ee=0$.


When $\lambda$ moves away from a special value then, in (\ref{eq:1.2}), $p$ remains at 0 while $q$ becomes small and nonzero.
We can then reduce (\ref{eq:1.2}) as in  Generic Case~1.1, as follows. The 3-jet of (\ref{eq:1.2}) becomes $(u_1, u_2, s_1^2+s_2^2u_2+s_2u_1^2+qs_2^3)$
with $q\ne 0$. Replacing $s_2$ by $ms_2+nu_2$ where $3qn+1=0$ and $qm^3=1$, and then removing terms in the third
component involving only $u_1, u_2$, reduces this to
\[ \left(u_1, u_2, s_1^2+\frac{1}{q^{1/3}}s_2u_1^2-\frac{1}{3q^{4/3}}s_2u_2^2+s_2^3\right).\]
The product of terms in front of $s_2u_1^2$ and $s_2u_2^2$ therefore has the sign of $-q$ and hence changes
as $q$ passes through 0. Furthermore it is not possible for both signs to be positive. We deduce the following.
\begin{cor}
Moving $\lambda$ through a special value $\lambda=g_3/(g_3\pm f_3)$ but keeping $\ee=0$
the type of equidistant always changes between Subcase~1.1.2 (negative definite) and Subcase~1.1.3 (indefinite)
as in Proposition~\ref{prop:def-indef}.  It is not possible to realize the positive definite Subcase~1.1.1.
\label{cor:special-unf}
\end{cor}
Figure~\ref{fig:special-clock} shows a typical way in which equidistants near to a special value evolve as
$\lambda$ and $\ee$ change.

\subsection{Some further details of Special Case 1.2}
We take $\lambda_0=\frac{g_3}{g_3+f_3}$ as a special value, assuming $f_{20}\ne 0, g_{20}\ne 0, f_3>0, g_3>0,
\gl_0 f_{20} + (1-\gl_0)g_{20} \ne 0$, i.e.\ $f_{20}g_3+g_{20}f_3 \ne 0$,  and also $R\ne 0$ (\ref{eq:R}) hold.
We write $\lambda=\lambda_0+\alpha$ for nearby values, and examine the full versal unfolding $\wH$ of $H$, as follows.

Thus the family of equidistants can be reduced to
\begin{equation}
 \wH(s_1,s_2,u_1,u_2,p,q)= \left(u_1, u_2, \pm s_1^2+s_2^2u_2+s_2u_1^2+ s_2^4+s_2^3u_1 + ps_2 + qs_2^3\right)=(u_1,u_2,\wh),
\label{eq:Htilde}
\end{equation}
say, where $p,q$ are unfolding parameters that are closely related to $\ee, \alpha$ respectively.

As an aid to understanding the equidistants for $(\ee,\alpha)$ close to $(0,0)$ we can calculate the loci in the $(p,q)$-plane at which the
structure of the singular set or the self-intersection set on the equidistant changes.

\begin{figure}
\begin{center}
\scalebox{1.4}{\includegraphics[width=1.4in]{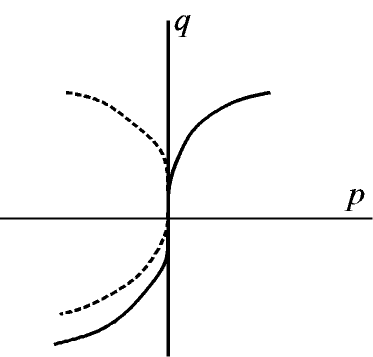}}
\end{center}
\vspace*{-0.2in}
\caption{\small Special Case 1.2. A schematic drawing of two curves in the $p,q$-plane at which the structure of the equidistant in the family~(\ref{eq:Htilde}) changes, either because the
cuspidal edge set changes (solid curve, together with the $q$-axis) or the self-intersection set changes (dashed curve). }
\label{fig:pq-plane}
\end{figure}

\begin{enumerate}
\item {\bf Singular set} \ For fixed $p,q$ the singular set is the image under $\wH$ of the set of points (using suffices for partial derivatives)
\[  (0,s_2,u_1,u_2) \mbox{ such that } \wh_{s_2}=\wh_{s_2s_2}=0.\]
Eliminating $u_2$, the equations reduce to
\[  u_1^2 - 3s_2^2u_1 + (p - 3s_2^2q-8s_2^3) = 0,\]
and the condition for this to have real solutions for $u_1$ is
\[ 9s_2^4 +32s_2^3 + 12qs_2^2 - 4p \ge 0.\]
We are therefore interested in finding the pairs $(p,q)$ for which there is a change in the number of real intervals in the set of $s_2$
satisfying this inequality.  This will occur when the discriminant with respect to $s_2$ vanishes, and that gives a locus of the form
\begin{equation}
p=0 \mbox{ or } p = \textstyle{\frac{1}{16}}\displaystyle q^3 + \textstyle{\frac{9}{1024}}\displaystyle q^4 + \ldots.
\label{eq:pqcuspedge}
\end{equation}
See Figure~\ref{fig:pq-plane}.

\item {\bf Self-intersection locus} \  Suppose $ (0, s_{21}, u_{1}, u_{2})$ and $(0, s_{22}, u_{1}, u_{2})$ are both
in the critical set of $\wH$ ($h_{s_1}=0$ gives $s_1=0$) and have the same image under $\wH$. Then with a little more trouble we can eliminate the $u$ variables and
obtain a condition in $s_{21}, s_{22}$ alone.  It is slightly more convenient to write $s_{21}=v_1+v_2, \ s_{22}=v_1-v_2$; then in fact we require
$ v_1(4v_1^3 +16v_1^2+8qv_1 + p+q^2) \ge 0$.  The number of $v_1$-intervals on which this holds will change when the discriminant with
respect to $v_1$ vanishes. One case here gives the same condition as (i) above, but we are concerned with
the remaining possibility: taking into account that $v_1, v_2$ must both have real solutions the locus in the $(p,q)$-plane is
\begin{equation}
p = -q^2, \ q\ge 0,
\label{eq:pqselfint}
\end{equation}
where of course the double root is $v_1=0$, that is $s_{22}=-s_{21}$. (The other potential double root when $p=-q^2$ leads to $q=2$ and is therefore
not relevant to a neighbourhood of the origin in the $(p,q)$-plane.)  See Figure~\ref{fig:pq-plane}.
\end{enumerate}

\section{Degenerate Case 2}\label{s:degen}
In this section we give some details of Degenerate Case~2, that is  $\gl f_{20} + (1-\gl)g_{20} = 0$. This gives a unique value of $\lambda$, namely $\displaystyle{\gl=\frac{g_{20}}{g_{20}-f_{20}}}$.
(If $f_{20}=g_{20}$ then, using $\gl f_{20} + (1-\gl)g_{20} = 0$, it follows that $f_{20}=g_{20}=0$, contrary to our assumptions.)  Thus whatever surfaces $M_0, N_0$ we start with
there will be an equidistant which falls into this case. It turns out to be a rich area for investigation; here we shall give some invariants which help to separate out the many subcases.
One of these invariants classifies the effect of changing $\lambda$ slightly from the degenerate value, while preserving the geometrical situation of two surfaces with parallel tangent
planes at parabolic points where the asymptotic directions are parallel, that is $\ee=0$ in (\ref{eq:M}), (\ref{eq:N}).  See Proposition~\ref{prop:case2near}.

\subsection{A normal form for Degenerate Case 2}
The 2-jet of $H_{0\gl}$ is now $(u_1,u_2, 2f_{20}s_1u_1)$.  Writing the third component as $u_1(s_1 + \mbox{h.o.t.}) + $ terms independent of $u_1$
and then using the bracketed expression to redefine $s_1$ we can eliminate $u_1$ from the higher terms.  Then replacing $s_2$ by
an expression of the form $s_2+au_2$ we can remove the degree 3 terms $s_1u_2^2$ and $s_1s_2u_2$. When this is done, the coefficient
of $s_2^2u_2$ becomes $3g_{030}f_{20}^2/g_{20}^2 \ne 0$ and the coefficient of $s_2u_2^2$ becomes $3f_{20}g_{030}(g_{20}-f_{20})/g_{20}^2\ne 0$.
We shall also assume that the coefficient of $s_1^3$ is nonzero to avoid further degeneration. We can now use scaling to reduce the 3-jet of
$H_{0\gl}$ to
\[ \left(u_1, \ u_2, \ s_1u_1+ s_1^3+s_2^2u_2+ s_2u_2^2+bs_1^2s_2+cs_1^2u_2+ds_1s_2^2+es_2^3\right),\]
for coefficients $b, c, d, e$.  The 4-jet can then be reduced by similar arguments, including scaling, to
\begin{eqnarray}
(u_1,u_2,h)&=&  \left(u_1, \ u_2, \ s_1u_1+ s_1^3+s_2^2u_2+ s_2u_2^2+bs_1^2s_2+cs_1^2u_2+ds_1s_2^2+es_2^3+ s_1^4 \right. \nonumber \\
&& \left. + (ps_2+qs_1^2) \right),
\label{eq:s1u1}
\end{eqnarray}
provided the cofficient of $s_1^4$ is nonzero: this and the 4-$\mathcal A$-determinacy of this 4-jet hold generically, by standard calculations.
The terms in brackets, $ps_1+qs_1^2$, represent an $\mathcal A$-versal unfolding of this germ.
We have not been able to reduce the number of coefficients $b,c,d,e$.  We shall work with (\ref{eq:s1u1}) as a `normal form'
and when appropriate interpret the coefficients in terms of the surfaces $M_0, N_0$.

The equidistant for $M_0, N_0$ and $\gl = g_{20}/(g_{20}-f_{20})$ is then locally diffeomorphic to the image under (\ref{eq:s1u1}) of the set $\{(s_1,s_2,u_1,u_2): h_{s_1}=h_{s_2}=0\}$.
Here, $h_{s_1}=0$ defines $u_1$ as a smooth function of the other three variables, while $h_{s_2}=0$ can be written
\begin{equation}
\frac{\partial h}{\partial s_2}=(s_2+u_2)^2+bs_1^2+2ds_1s_2+(3e-1)s_2^2=(s_2+u_2)^2-T(s_1,s_2)=0,
\label{eq:dhds2}
\end{equation}
say where $T$ is a quadratic form in $s_1, s_2$ which we shall assume to be nondegenerate, that is $d^2-b(3e-1)\ne 0$.

\subsection{Plotting the equidistants}\label{ss:param}

It is also useful to rewrite the equation of the quadric cone $C$, given by $h_{s_2}=0$, where $p=q=0$ in (\ref{eq:s1u1}), and provided $b \ne 0$, as
\begin{equation}
C: \ \ (s_2+u_2)^2 + b\left(s_1+\frac{d}{b}s_2\right)^2 + \left(\frac{3be-b-d^2}{b}\right)s_2^2=0.
\label{eq:cone}
\end{equation}
Note that this is a single point at the origin if and only if all coefficients are $>0$ (since the first one is $>0$),
that is
\[ b> 0, \ d^2+b-3be<0;\]
compare Proposition~\ref{prop:case2near}.

 The equidistant (for $p=q=0$)
is the image of $C$ under the map $\RR^3\to\RR^3$ given by
\[ (s_1,s_2,u_2)\mapsto (u_1, u_2, \overline{h}(s_1,s_2,u_2))\]
where on the right-hand side $u_1$ is expressed in terms of $s_1,s_2,u_2$
using $h_{s_1}=0$ and this is substituted into $h$, giving the function $\overline{h}$.

 We can find a `good' parametrization of the equidistant by using coordinates $(x_1,x_2,s_2)$ and writing (\ref{eq:cone}) as
\[ x_1^2 + bx_2^2 + ks_2^2, \mbox{ where } k = \textstyle{\frac{3be-b-d^2}{b}}, \ x_1=s_2+u_2, \ x_2=s_1+\textstyle{\frac{d}{b}}\displaystyle s_2.\]
Thus the substitution to use in $\overline{h}$ is $u_2=x_1-s_2, \ s_1=x_2-\textstyle{\frac{d}{b}}\displaystyle s_2$. The equidistant
is then plotted as follows.

\begin{enumerate}
\item
If $b>0$ and $C$ is not a single point then  $k<0$ (i.e.\ $d^2+b-3be>0$) and we write
\[ x_1^2 + bx_2^2 = (-k)s_2^2,\]
so that for any $(x_1, x_2)\ne(0,0)$ we have two distinct values for $s_2$: there is no restriction on the values of $x_1,x_2$. We use
$x_1,x_2$ as parameters and the two `halves' of $C$ are given by the two values of $s_2$.
\item
If $b<0, k>0$ (i.e.\ $ d^2+b-3be>0$)  then we similarly write $x_1^2 +k s_2^2 = (-b)x_2^2,$
so that for any $(x_1, s_2)\ne(0,0)$ we have two distinct values for $x_2$.
Here $x_1,s_2$ are used as parameters.
\item
Finally if $b<0, k<0$ (i.e.\ $d^2+b-3be<0$) then we write $x_1^2 = (-b)x_2^2 + (-k)s_2^2$ and for any $(x_2, s_2)\ne(0,0)$ we have two distinct values for $x_1$.
Here $x_2,s_2$ are used as parameters.
\end{enumerate}

For values of $(p,q)$ other than $(0,0)$ the equation of $C$ acquires an extra term $-p$ on the right-hand side,
thus creating a hyperboloid of one or two sheets (or an ellipsoid when $C$ is a single point). In fact the
hyperboloid has one sheet when $bkp >0$, that is $(d^2+b-3be)p<0)$, and two sheets when $bkp<0$, that is $(d^2+b-3be)p>0)$. In the two-sheet situation
the same method as above plots the equidistant, without restrictions on the values of the parameters.
In the one-sheet situation the points in the parameter plane lie outside an ellipse, the `waist' of the hyperboloid.  This ellipse is given in the three situations above
by $ x_1^2 + bx_2^2 =-p, \ x_1^2 +k s_2^2=-p$ and $(-b)x_2^2 + (-k)s_2^2=p$ respectively.
 In the situation where $C$ is
a single point, and $p<0$, the points in the parameter plane lie inside an ellipse. In all situations, $q$ does not
affect the hyperboloid or ellipsoid, but of course its value affects the function $\overline{h}$.

\subsection{Nearby non-special values of $\lambda$}\label{ss:case2nearby}
Here, we examine the effect of adding in the term $qs_1^2$ in (\ref{eq:s1u1}).  This represents changing $\lambda$ from the value $g_{20}/(g_{20}-f_{20})$ to a
nearby value, which will be of the type considered in Generic Case~1.1, provided the coefficient $e$ of $s_2^3$ in (\ref{eq:s1u1}) is nonzero, and to avoid further
degeneracy we shall  assume this to be true.
We determine here, in terms of $b,c,d,e$, which subcase of Proposition~\ref{prop:def-indef} is obtained, and then refer this back to the surfaces $M_0,N_0$. (The subcase does
not depend on the sign of $q$ in the added term $qs_1^2$.) To do  this we reduce (\ref{eq:s1u1}), with $p=0$ but with $qs_1^2$ present, to the normal form
found above for Generic Case~1.1,  by making the `left' and `right' changes of coordinates as sketched above. We can restrict attention for this to the terms of (\ref{eq:s1u1}) of degree
$\le 3$ since the Generic Case~1.1 germ is 3-$\mathcal A$-determined. Thus we start by redefining $s_1$ (`completing the square') to change the degree 2 terms to $s_1^2$,
remove the terms in $u_1, u_2$ only, remove the remaining terms besides $s_1^2$ that are divisible by $s_1$ and then redefine $s_2$ by adding suitable multiples
of $u_1$ and $u_2$. The result of this is to reduce the 3-jet of  (\ref{eq:s1u1}) by $\mathcal A$-equivalence to the form
\[ \left( u_1, \ u_1, \ qs_1^2 + e s_2^3 + \frac{s_2}{12eq^2}\left( (3be-d^2)u_1^2+4qdu_1u_2+4q^2(3e-1)u_2^2 \right)\right). \]
The discriminant of the quadratic form in $u_1, u_2$ is $(d^2+b-3be)/3eq^2$, so this form is definite if and only if $e(b+d^2-3be)<0$.  Scaling so that the terms in
$s_1^2, s_2^3$ have coefficients equal to 1 multiplies the quadratic form in $u_1, u_2$ by $(q^2e)^{-1/3}$, and from this we deduce the following, where (i) and (ii) are
derived by direct calculations from the parametrizations of $M_0$ and $N_0$.
\begin{prop}\label{prop:case2near}
The normal form (\ref{eq:s1u1}) for Degenerate Case~2, with $p=0$ but $q$ nonzero and small, corresponding to a small change in $\lambda$,
gives the following subcases of Generic Case~1.1 (general $\gl$): \\
Subcase 1.1.1 (positive definite, $++$): $e>\frac{1}{3}$ and $d^2+b-3be < 0$,\\
Subcase 1.1.2 (negative definite, $- -$): $e<\frac{1}{3}$ and $e(d^2+b-3be) < 0$,\\
Subcase 1.1.3 (indefinite, $+-$): $e(d^2+b-3be > 0$.

\medskip\noindent
In terms of the surfaces $M_0, N_0$, \\
{\rm (i)} \ When $f_{030}g_{030}> 0$, so $f_{030}=f_3^2, \ g_{030}=g_3^2$, $e<\frac{1}{3}$ and has the sign of $f_{20}g_3^2-g_{20}f_3^2$ while $d^2+b-3be$ has the sign of $-R$ as in (\ref{eq:R}).\\
{\rm (ii)} \ When $f_{030}g_{030}< 0$, so $f_{030}=f_3^2, \ g_{030}=-g_3^2$, $e>\frac{1}{3}$ and $d^2+b-3be$ has the sign of $R$.
\end{prop}

\subsection{Invariants distinguishing subcases of Degenerate Case~2}\label{ss:invariants}
We shall use the following:\\
\begin{enumerate}
\item The number of cuspidal edges on the equidistant for $p=q=0$, which can be 0, 2 or 4 (see below);
\item The number of self-intersection curves on the equidistant for $p=q=0$, which can be 0, 1, 2 or 3 (see \S\ref{ss:SIdegen});
\item The subcase of Generic Case~1.1 given in Proposition~\ref{prop:case2near} which is obtained by changing $\lambda$ slightly.
\end{enumerate}
This might give $3 \times 4 \times 3 = 36$ subcases but fortunately many of these combinations cannot be realized. We shall give
values of $b,c,d,e$ realizing  of all possible subcases in \S\ref{ss:examples}, Table~\ref{table1} below.

For given values of these invariants, the interval in which $e$ lies, either $e<0$ or $0<e<\frac{1}{3}$ or $e>\frac{1}{3}$ could in principle affect the equidistant but so
far as we are aware the basic geometrical structure---the qualitative nature of the equidistant---is not affected.

\medskip

The number of cuspidal edges, that is 1-dimensionial singular sets, on the equidistant, can be calcuated as follows.
We can regard $h_{s_2}=0$, as in \S\ref{ss:param} above, as the equation of a quadric cone $C$ in $\mathbb R^3$ with coordinates $(s_1, s_2, u_2)$.
The quadric cone $C$ is nondegenerate since $T$ in (\ref{eq:dhds2}) is a nondegenerate quadratic form,
and consists of the origin alone if and only if $T$ is negative definite (that is, $d^2<b(3e-1)$ and $b>0$), otherwise it is a real cone, or equivalently a real nonsingular conic in $\RR P^2$.

When $T$ is not negative definite, the equidistant therefore has two `branches', which are the images of the two halves of the cone; these branches may
intersect (apart from at the origin) and will generally  themselves be singular.  Writing the equation of $C$ more briefly as $\g(s_1,s_2,u_2)=0$,
the singular set of the equidistant is the image of certain curves on $C$, given by
the additional equation \[ \overline{h}_{s_1}\g_{s_2}-\overline{h}_{s_2}\g_{s_1}=0.\]
(This can be written in terms of $h$ itself as $h_{s_1s_1}h_{s_2s_2}-h_{s_1s_2}^2=0$.)
The lowest terms of the left hand side are of degree 2 in $s_1, s_2, u_2$ and therefore give another conic $C_2$
in $\RR P^2$. The equation of $C_2$
is in fact
\[ (b^2-3d)s_1^2+(bd-9e)s_1s_2-(cd+3)s_1u_2+(d^2-3be)s_2^2-(3ce+b)s_2u_2-cu_2^2 = 0.\]
 This meets the nonsingular conic $\g = 0$ in 0, 2 or 4 real points. (The conic $C_2$ cannot in fact be a single point: examination of the
matrix of the above quadratic form in variables $s_1,s_2,u_2$ defining $C_2$ shows that its determinant is always $\le 0$ so the quadratic form cannot be positive
definite, and negative definiteness is also ruled out by examining the signs of the other leading minors. The leading $1\times 1$ minor cannot be $<0$ at the
same time as the leading $2\times 2$ minor is $>0$.)
There are therefore 0, 2 or 4 curves through the origin on $C$ whose images are the singular points,  the cuspidal edges, of the equidistant. These cuspidal edges
pass through the origin, lying on both `sheets' of the equidistant.

The number of cuspidal edges can be calculated for example by substituting $(s_1,s_2,u_2)=
(mt,nt,t)$ in the  equations of  $C$ and $C_2$, taking out the factor $t^2$ and finding the common solutions of the two resulting quadratic equations in $m,n$.
Eliminating one of $m,n$ gives a degree 4 equation in the other and there are standard algebraic techniques for computing the number of real solutions of
a quartic equation---or for given $(b,c,d,e)$ we can solve numerically.  The results for the Classes~I-X are given in Table~\ref{table1} below.

\subsection{Self-intersections of the equidistant in Degenerate Case 2}\label{ss:SIdegen}

We start with the normal form (\ref{eq:s1u1}) in \S\ref{s:degen}, namely
\[
(u_1,u_2,h) = \]
\[ \left(u_1, \ u_2, \ s_1u_1+ s_1^3+s_2^2u_2+ s_2u_2^2+bs_1^2s_2+cs_1^2u_2+ds_1s_2^2+es_2^3+ s_1^4 +ps_2+qs_1^2 \right),
\]
subject to the critical set conditions $h_{s_1}=h_{s_2}=0$.  We include the unfolding terms $ps_2+qs_1^2$ though we are particularly interested
in the self-intersections for $p=q=0$.  We can immediately solve $h_{s_1}=0$  for $u_1$:
\[ u_1= -2bs_1s_2-2cs_1u_2-ds_2^2-3s_1^2-4s_1^3-2qs_1,\]
so that the equations which state that two domain points $(s_1,s_2,u_1,u_2)$ and say $(t_1,t_2,u_1,u_2)$ have the same image
take the following form.

\noindent
(SI1): the above formula for $u_1$ gives the same answer for both domain points;\\
(SI2): the formula for $h$ above gives the same answer for both domain points;  \\
(SI3): $h_{s_2}(s_1,s_2,u_1,u_2)=0$; and\\
(SI4): $h_{t_2}(t_1,t_2,u_1,u_2)=0.$

It is convenient to make the substitution $s_1=x_1+y_1, t_1 = x_1-y_1, s_2=x_2+y_2, t_2=x_2-y_2$, so that the
`trivial solution' $s_1=t_1, s_2=t_2$ becomes $y_1=y_2=0$. Furthermore replacing $y_1$ by $-y_1$ and $y_2$ by $-y_2$ interchanges
$(s_1,s_2)$ and $(t_1,t_2)$, that is interchanges the two domain points $(s_1,s_2,u_1,u_2)$ and $(t_1,t_2,u_1,u_2)$ with the same image
in $\RR^3$ under the normal form map  (\ref{eq:s1u1}) . With this substitution the equations become say (SI1$'$), etc., and we
use (SI3$'$)-(SI4$'$) to solve for $u_2$:
\[ u_2=-\frac{bx_1y_1+dx_1y_2+dx_2y_1+3ex_2y_2}{y_2}, \]
where the denominator $y_2$ is harmless since it is easy to check that if $y_2=0$ then the other equations imply
that $y_1=0$ too.  {\em Note that this expression does not involve} $p,q$.

\begin{figure}[!ht]
\begin{center}
\scalebox{1.4}{\includegraphics[width=1.3in]{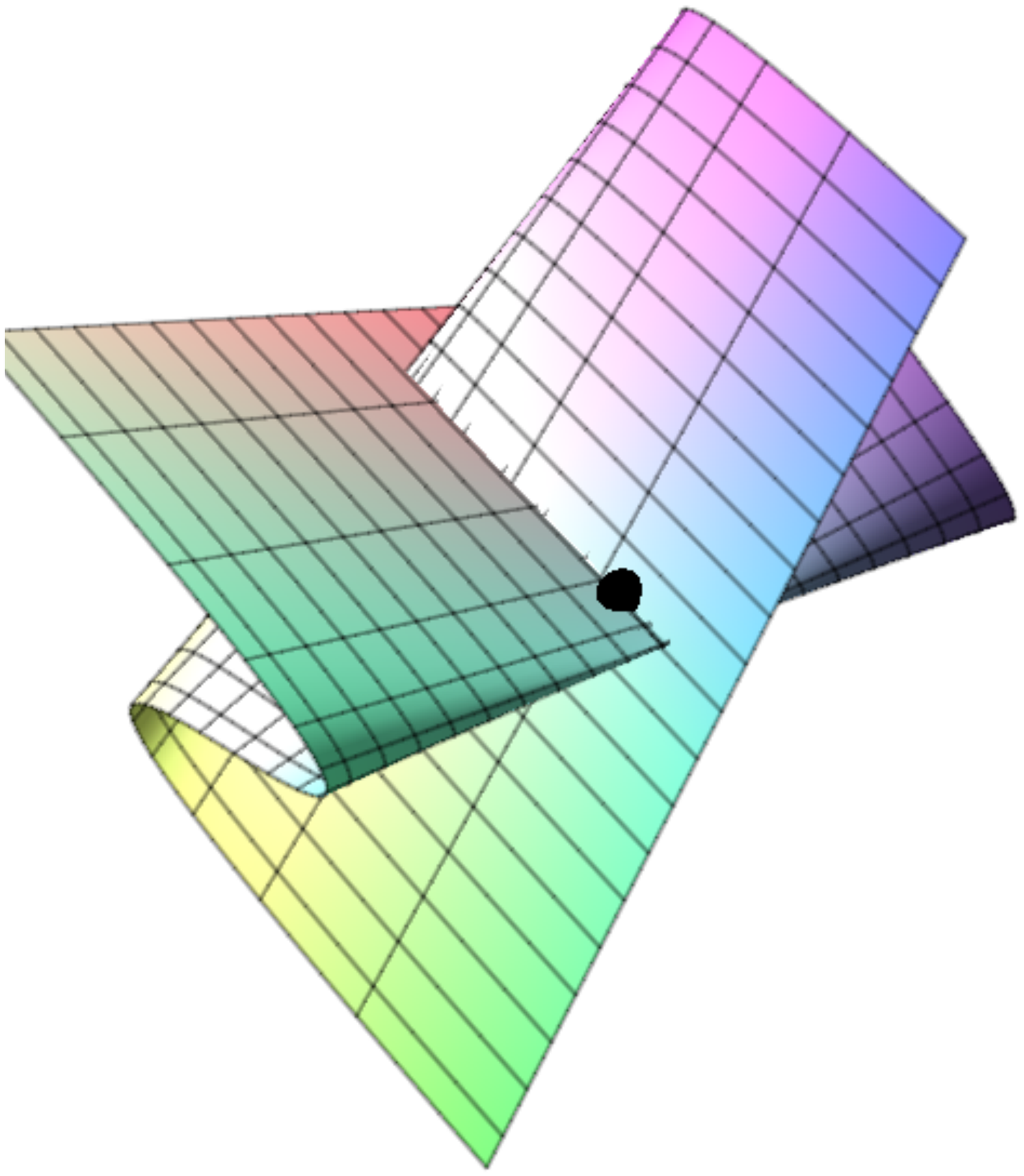}}
\scalebox{1.4}{\includegraphics[width=1.3in]{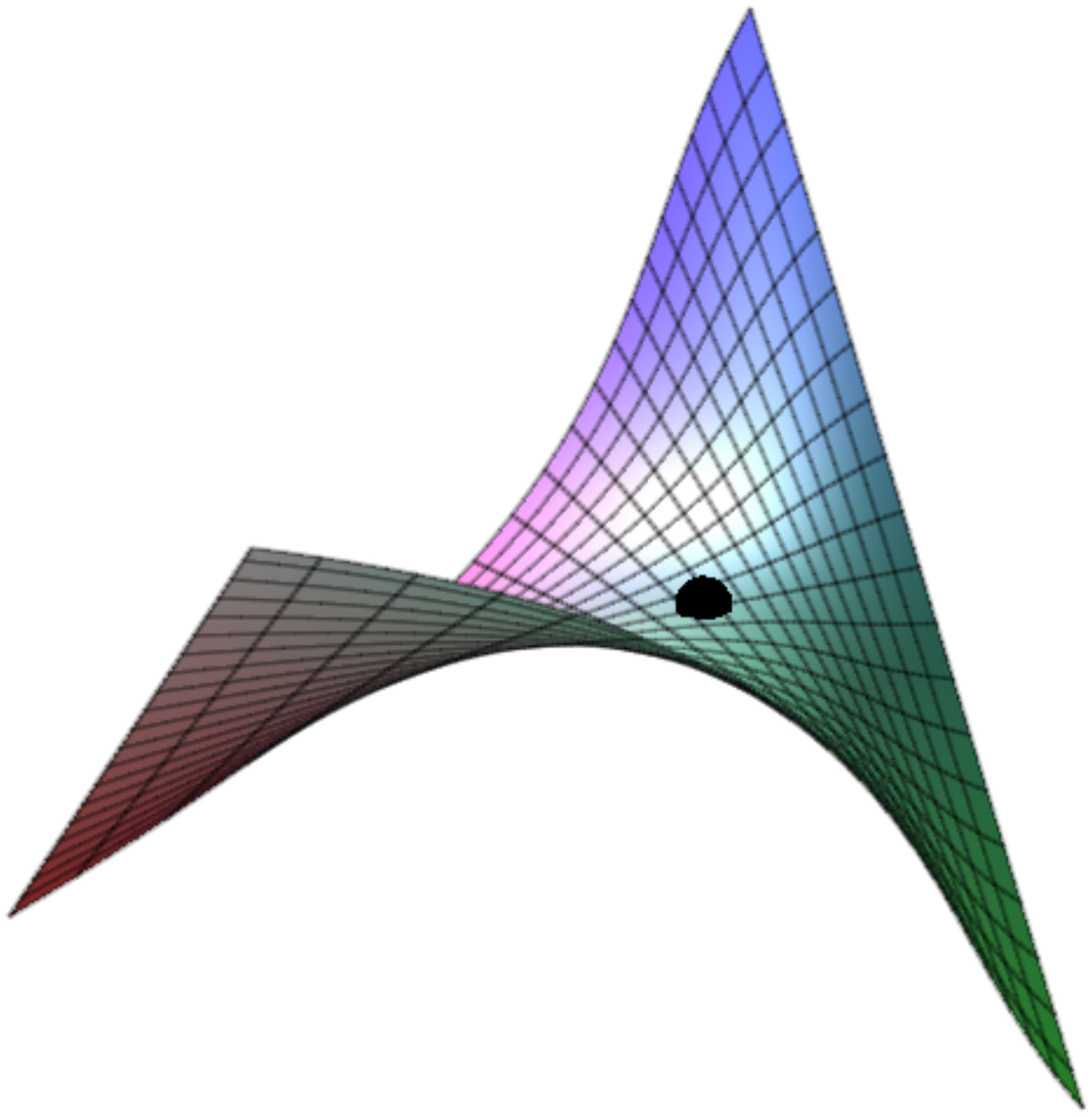}}
\end{center}
\vspace*{-1cm}
\caption{\small  The surface given by (\ref{eq:SI5}) or (\ref{eq:SI5param}), for (left) $b=8, c=-4, d=-3, e=-1$, with three smooth sheets
through the origin, which is marked by a black dot; (right) $b=-8, c=4, d=-3, e=-1$, with one smooth sheet. (See Proposition~\ref{prop:sheets}.)
These are respectively
Class~III and Class~IX in Table~\ref{table1} below. Note that in the first of these there are nevertheless only two self-intersection curves of the equidistant
for $p=q=0$ , using the
criterion of Proposition~\ref{prop:self-int-number}.  In fact the picture for Class II is very similar to the left-hand figure, but  there is only
one self-intersection curve of the equidistant for $p=q=0$.}
\label{fig:SI5ab}
\end{figure}

We can solve (SI1$'$) for $x_2$:
\[ x_2 = \frac{bcx_1y_1^2+cdx_1y_1y_2-bx_1y_2^2-6x_1^2y_1y_2-2y_1^3y_2-3x_1y_1y_2-qy_1y_2}
{-cdy_1^2-3cey_1y_2+by_1y_2+dy_2^2}.\]
This time we may need to investigate the vanishing of the denominator, but
assuming the denominator is nonzero and substituting for $x_2$ we find that
the equation (SI2$'$)-$y_2$((SI3$'$)+(SI4$'$))  reduces to
\begin{equation}
 \mbox{SI5}: by_1^2y_2+dy_1y_2^2+ey_2^3+4x_1y_1^3+y_1^3=0.
 \label{eq:SI5}
 \end{equation}
This is to be treated as the equation of a surface in 3-space $(x_1,y_1,y_2)$ which contains the
$x_1$-axis, since $(x_1,0,0)$ is always a solution.  The surface will have a certain number of `sheets'
passing through the origin, equal to the number of values of $k$ which make the first coordinate zero in the following parametriztion of SI5 by $k$ and $y_1$.
\begin{equation}
\left( -\frac{ek^3+dk^2+bk+1}{4}, \ y_1, \ ky_1\right).
\label{eq:SI5param}
\end{equation}
If $y_1=0$ in (\ref{eq:SI5}), then $y_2=0$ and $x_1$ is arbitrary; and indeed, being cubic in $k$, (\ref{eq:SI5param}) gives all
points $(x_1,0,0)$, possibly for more than one (real) $k$.  If $y_1\ne 0$ then we solve (\ref{eq:SI5}) for $x_1$ and writing $y_2=ky_1$ produces the given
value $-\frac{1}{4}(ek^3+dk^2+bk+1)$ for $x_1$.  Conversely, every point (\ref{eq:SI5param}) satisfies (\ref{eq:SI5})
by substitution.  Hence (\ref{eq:SI5param}) parametrizes the complete surface (\ref{eq:SI5}).  Two examples are shown in
Figure~\ref{fig:SI5ab}.

Note that the surface (\ref{eq:SI5}) and the parametrization (\ref{eq:SI5param}) are independent of the unfolding parameters $p,q$.

\begin{prop}
The number of smooth real sheets of the surface (\ref{eq:SI5}) through the origin in $(x_1,y_1,y_2)$-space is 1 or 3 according as
\[ 27e^2 + 2b(2b^2-9d)e +d^2(4d-b^2)> 0 \mbox{ or } <0 \mbox{ respectively}.\]
This number is therefore the maximum number of self-intersection branches of the equidistant, for any $p,q$.
If $b^2<3d$ then the displayed expression is $>0$ for all values of $e$.
\label{prop:sheets}
\end{prop}
{\em Proof} \ This is a matter of calculating the discriminant of the cubic polynomial $ek^3+dk^2+bk+1$ in $k$, and the discriminant
$16(b^2-3d)^3$ of the displayed
quadratic polynomial in $e$.  The sheets will be smooth provided the cubic in $k$ has no repeated root, that is
provided the discriminant is nonzero. \hfill$\Box$

\begin{rem}\label{rem:D4-contact}
{\rm
In \S\ref{ss:contact} we noted that, in the current Degenerate Case~2, the sign of a certain polynomial in the coefficients of the two surfaces $M_0, N_0$ determines
whether the `scaled contact map' has type $D_4^+$ or $D_4^-$. By reducing to normal form as in \S\ref{s:normal} we can re-express this polynomial in terms
of the coefficients $b,c,d,e$ of the normal form.  When this is done, we find that the condition for one (resp.\ three) sheets as in the above proposition coincides
with the condition for $D_4^+$ (resp.\ $D_4^-$) in the scaled contact map.  We do not know the full significance of this fact.
}
\end{rem}

\medskip

Substituting  $x_1=-\frac{1}{4}(ek^3+dk^2+bk+1)$ and $y_2=ky_1$ in one of the conditions on $x_1, y_1, y_2$
not fully  used yet (for example, SI2$'$) we obtain a single equation
in $y_1, k$ (involving now $p$ and $q$) which determines the branches of the self-intersection set of the equidistant. We are interested in values of $k$ close
to a zero $k_0$ of the polynomial $ek^3+dk^2+bk+1$, so we now write $k=k_0+z$ say where $z$, as well as $y_1, p, q$, will be small. Since $k_0$ satisfies a cubic
equation we can express $k_0^3$ in terms of $k_0$ and $k_0^2$, namely as $k_0^3 = (-dk_0^2 - bk_0 - 1)/e$, and therefore all higher powers of $k_0$ can
be expressed in terms of $k_0, k_0^2$ as well.

\begin{defs}{\rm
For a chosen value of  $k_0$, the polynomial in $y_1, z, p, q$ just formed, the zero set of which determines
the solutions to (SI1)-(SI4) or their equivalents (SI1$'$)-(SI4$'$), and hence determines the points
corresponding to self-intersections of the equidistant,
 will be called $L(k_0)$. In the special case $p=q=0$, we shall write $L_0(k_0)$ for the polynomial in $y_1$ and $z$.
 }
\end{defs}
We deduce the following; the statements 2-5 are easily checked by direct calculation.

\begin{prop}
\begin{enumerate}
 \item For each real root $k_0$ of $ek^3+dk^2+bk+1=0$ one smooth sheet of the surface (\ref{eq:SI5}) is parametrized by $(y_1,z)$
and the points which correspond to self-intersections on the equidistant for any $p,q$ are given by the additional equation $L(k_0)=0$.

\item The polynomials $L(k_0)$ and $L_0(k_0)$ contain only the powers $y_1^2$ and $y_1^4$ of $y_1$.  For any $p,q$ the zero-set of $L(k_0)$ is
symmetric about the $y_1$-axis in the $(y_1,z)$-plane.

\item
The other variable $z$ occurs to powers $\le 14$ in $L(k_0)$. The coefficient of $z^{14}$ is in fact $27e^5(3e-1)$ which will not be zero
since $e=0, \frac{1}{3}$ are excluded values.

\item
The linear part  of $L(k_0)$ has the form constant $\times p$. The nonzero quadratic terms  are in $y_1^2, z^2, zp, zq$ and $q^2$.

\item
The 2-jet of $L_0(k_0)$ has the form $c_0 y_1^2 + c_2 z^2$.
\end{enumerate}
\label{prop:L}
\end{prop}
\medskip

The last statement above implies that, for $p=q=0$, a given sheet of the surface (\ref{eq:SI5}), that is a given value of $k_0$,
will correspond to a branch of the self-intersection set of the equidistant if and only of $c_0,  c_2$ have opposite signs. When $c_0c_2>0$ there is only
an isolated point at $y_1=z=0$. When
$c_0c_2<0$ the two real branches of the set $L_0(k_0)=0$ (forming a crossing at the origin $y_1=z=0$) will give only one branch of the self-intersection set because, as noted above,
replacing $y_1$ by $-y_1$, and hence $y_2=ky_1$ by  $-y_2=k(-y_1)$, merely interchanges the domain points contributing to the self-intersection.

Each of $c_0, c_2$ is quadratic in $k_0$; multiplying them gives an expression of degree 4 which can be reduced to degree 2 again using the equation $ek^3+dk^2+bk+1=0$.  Writing the
resulting quadratic expression as $N=N_0(b,c,d,e)+N_1(b,c,d,e)k_0+N_2(b,c,d,e)k_0^2$ we have the following, which is used to determine the number of self-intersection branches of
the equidistant in the ten classes of Table~\ref{table1}.
\begin{prop}
The number of real branches of the self-intersection set of the equidistant for $p=q=0$ is the number of solutions $k=k_0$ of $ek^3+dk^2+bk+1=0$ at which the quadratic $N$ is $<0$.
\label{prop:self-int-number}
\end{prop}

As $(p,q)$ moves away from $(0,0)$ we can still trace the zero set of $L(y_0)$ in the $(y_1,z)$-plane. An isolated point may disappear or open into a symmetric loop, which represents
a self-intersection of the equidistant having two endpoints, if the loop crosses the $y_1$-axis, and a closed self-intersection curve if it does not.
 A crossing will become a `hyperbola'; if it crosses the $y_1$-axis then the
corresponding self-intersection curve will have two endpoints and if not then it will be an unbroken arc. This is illustrated in the next section.

\subsection{Examples}\label{ss:examples}
Considering different realizable values of the three invariants in~\S\ref{ss:invariants}, we have the  ten classes of equidistant given in Table~\ref{table1}.  It is also possible in
some of these classes to allow  values of $e$ in  different ranges $e<0, \ 0<e<\frac{1}{3}, \ e >\frac{1}{3}$ but this does not appear to affect the
equidistant in any qualitative way. We can compute the curves in the $(p,q)$-plane alomg which the cusp edges or the self-intersection curves on the equidistant underfgo a
qualitative change.  (The ten cases of the table in fact have ten distinct configurations of these curves.)

\begin{table}[!ht]
\begin{center}
\begin{tabular}{|c|c|c|c|r|r|r|r|}\hline
Class & Cusp edges & self-int &Subcase &$b$&$c$&$d$&$e$ \\
&&&  (Prop.~\ref{prop:case2near}) &&&& \\ \hline\hline
I & 0 & 0 & $++$ & 8 & 4 & $-3$ & 1 \\ \hline
II & 0 & 1 & $+-$ & 8 & $-4$ & $-3$ & $\frac{1}{6}$ \\ \hline
III & 0 & 2 & $- -$ & 8 & $-4$ & $-3$ & $-1$ \\ \hline
IV & 2 & 0 &$+-$ & $-13$ & 6 & $-3$ & $-5$ \\ \hline
V & 2 & 1 & $- -$ & 1&2&3& $-1$ \\ \hline
VI & 2&2& $+-$ & 8 & 4 & $-3$ & $\frac{1}{6}$ \\ \hline
VII &2&3&$- -$ & $-13$ & $-6$ &  1 & $\frac{1}{6}$ \\ \hline
VIII&2&3&$+-$&$-8$&4&1&$\frac{1}{6}$ \\ \hline
IX & 4 &1&$+-$&$-8$&4&$-3$&$-1$ \\ \hline
X& 4&3&$+-$&$-8$&6&$-3$&10 \\ \hline
\end{tabular}
\end{center}
\caption{\small Ten distinct classes of Case 2, giving all possible realizations of the three invariants of \S\ref{ss:invariants}, and examples of
values of $b,c,d,e$ which realize these invariants. The fourth column refers to the `non-special' type
which results from changing $\lambda$ slightly from the degenerate value.}
\label{table1}
\end{table}

\begin{figure}[!ht]
\begin{center}
\scalebox{1.4}{\includegraphics[width=3in]{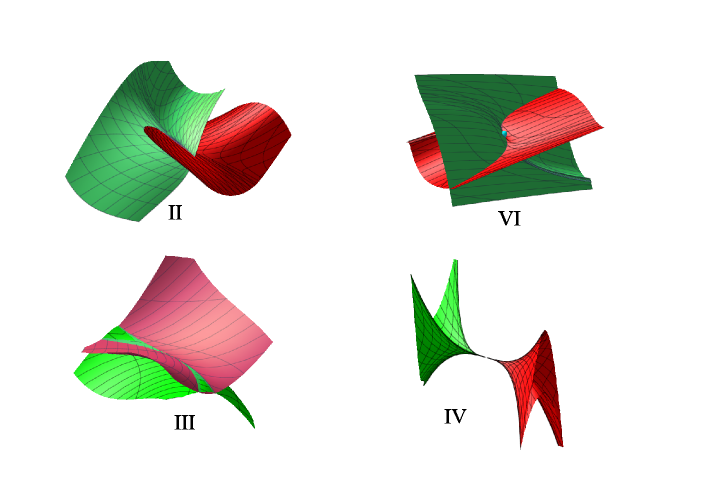}}
\end{center}
\vspace*{-1cm}
\caption{\small Cases II, III, IV and VI from Table~\ref{table1}, for $p=q=0$. The origin is marked for Case VI, where there are
two very arrow swallowtails passing through the origin, contributing two cusp edges and one self-intersection,  and the other
self-intersection is visible where the sheets pass through one-another. }
\label{fig:II-III-IV-VI}
\end{figure}

\begin{figure}[!ht]
\begin{center}
\scalebox{1.4}{\includegraphics[width=3in]{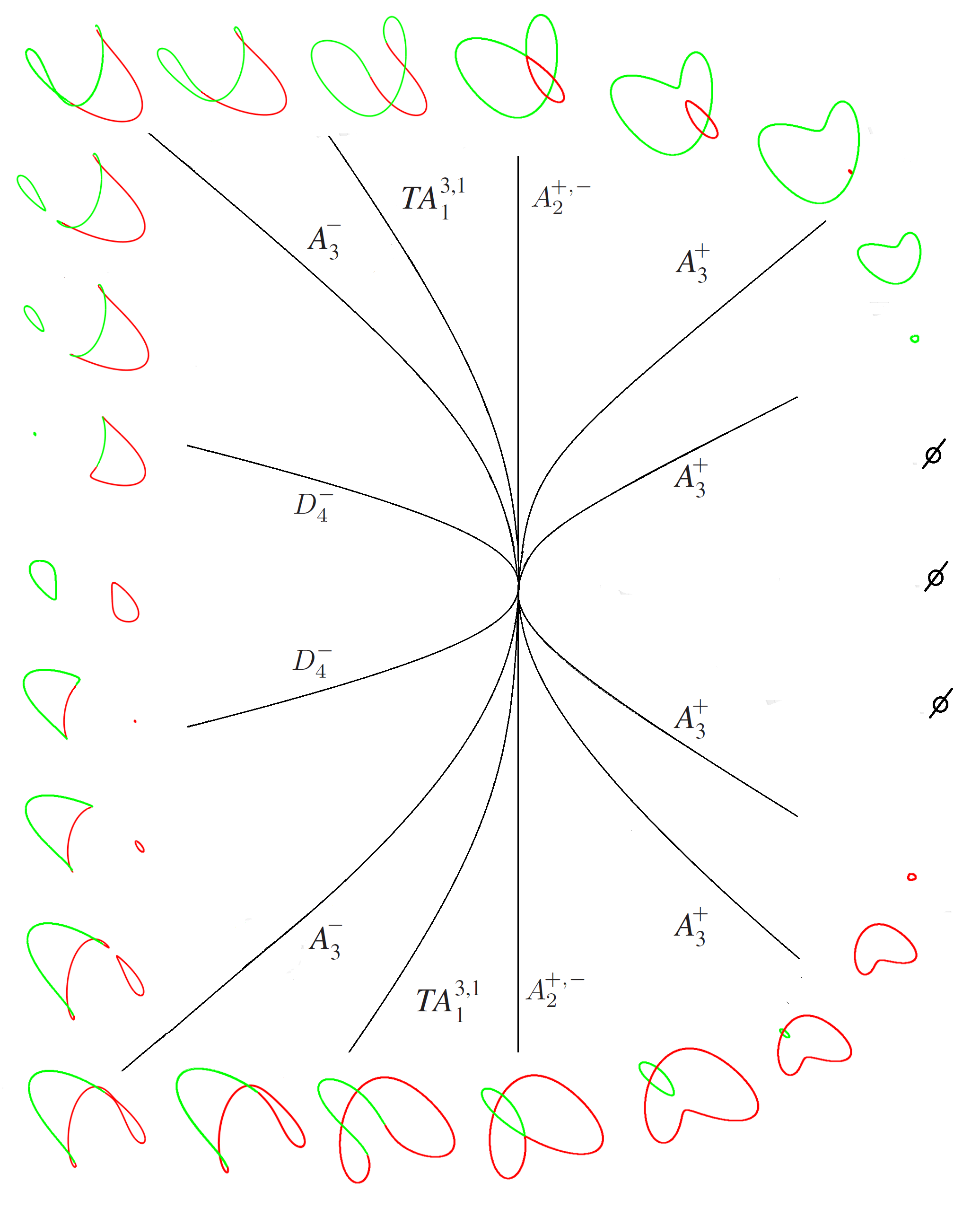}}
\end{center}
\vspace*{-1cm}
\caption{\small Pre-images of the cuspidal edges on the equidistants in Class II of Table~\ref{table1} for unfolding parameters $(p,q)$ making a circuit of the origin.
The colours correspond to either the two parts of a hyperboloid of two sheets as in \S\ref{ss:param} or to the two parts into which a hyperboloid of one sheet is
cut by the plane through the `waist'.  For the labelling of transitions, see Notation~\ref{not}.  }
\label{fig:Case9-clock5}
\end{figure}

\begin{figure}[!ht]
\begin{center}
\scalebox{1.4}{\includegraphics[width=4in]{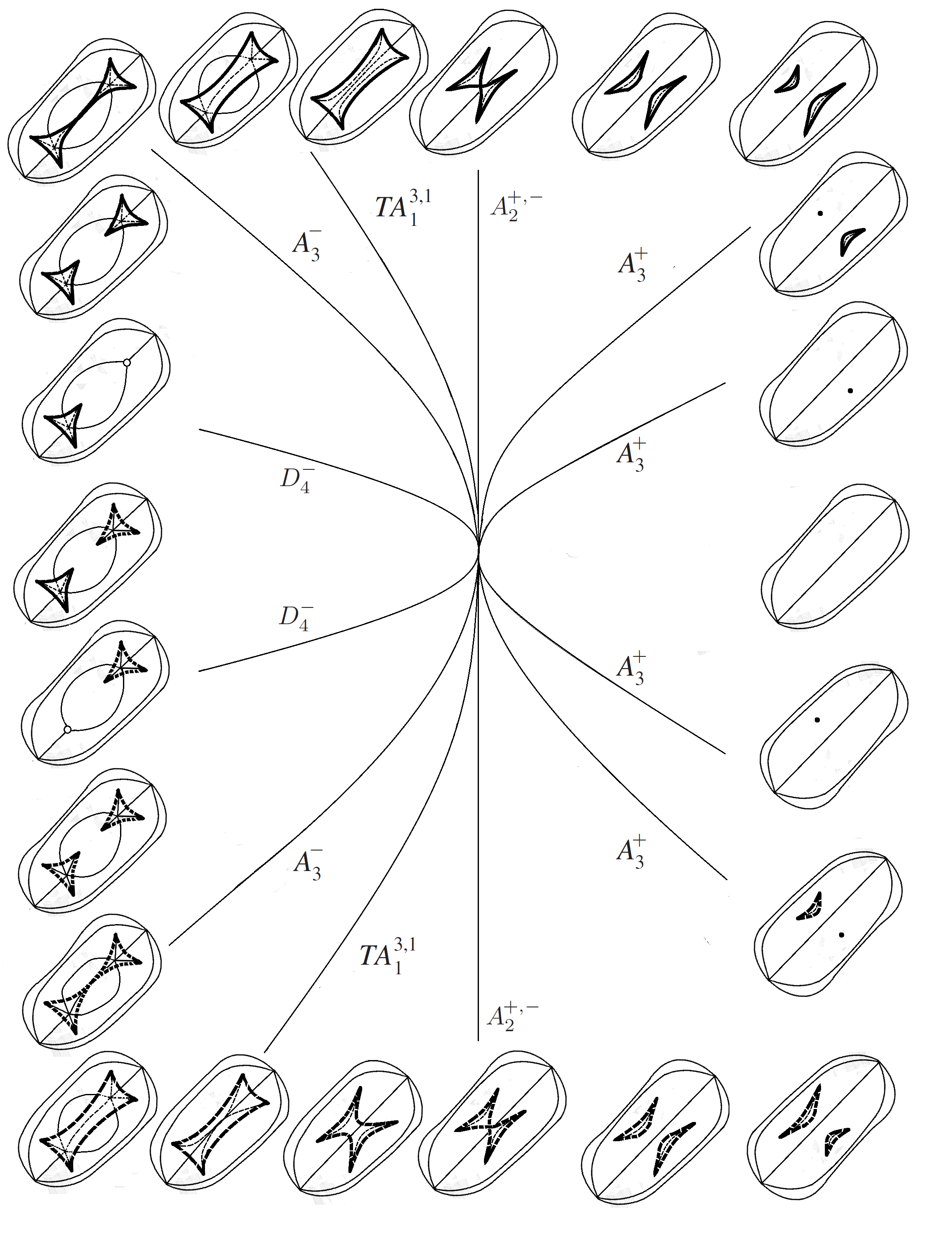}}
\end{center}
\vspace*{-1cm}
\caption{\small Schematic diagram of the equidistants for Class II of Table~\ref{table1}, with the unfolding parameters $(p,q)$ making a
circuit of the origin.  The figure shows cuspidal edges (thick lines) and self-intersections (thin lines) with
solid and dashed curves indicating visibility from one direction. For the labelling, see Notation~\ref{not}. }
\label{fig:Case9-WF-clock4.1}
\end{figure}

We  shall now give more detail on Case~II of the table, showing how the cuspidal edges and self-intersections of the equidistant evolve as $(p,q)$ in (\ref{eq:s1u1}) makes a
circuit of the origin.  Figure~\ref{fig:Case9-clock5} shows the transformations in the cuspidal edge as $(p,q)$ moves in such a circuit and
Figure~\ref{fig:Case9-WF-clock4.1} gives schematic diagrams of the corresponding
equidistants, indicating their self-intersections and cusp edges. We use the following labelling on these figures to indicate transitions (perestroikas) in the  structure
of the equidistant.

\begin{notation}\label{not}\verb++
{\rm
\noindent
$A_2^{++}, A_2^{- -}, A_2^{+ -}$ refer to Subcases 1.1.1, 1.1.2 and 1.1.3, as in Proposition~\ref{prop:def-indef}.  The corresponding transitions have also been desscribed to
as `Zeldovich's pancakes' or `flying saucer', the `hyperbolic transformation of an edge', and `the death of a compact component of an edge', respectively. See also \cite{Gory,Gory-S}.

\smallskip\noindent
$A_3^+, A_3^-$ refer to the `swallowtail-lips' and `swallowtail-beaks'  singularity respectively.

\smallskip\noindent
$D_4^-$ refers to the `pyramid' singularity (and  $D_4^+$ would similarly be the `purse' singularity).

\smallskip\noindent
$TA_1^{3,1}$, called such in \cite{Gory-S,Sul} (see also \cite{Gory}) refers to the situation where three smooth sheets of the
equidistant are pairwise transversal to each other, but the curve
of intersection of any two of them is tangent to the third sheet at the moment of bifurcation.
}
\end{notation}

\section{Conclusion and further work}\label{s:conc}

There have been many recent studies of singularities of (affine) equidistants of surfaces.  For a single equidistant of
a fixed surface, the generic singularities are $A_1, A_2, A_3$ (see for example \cite{GZ,domitrz2});
for a fixed surface, but allowing the ratio $\lambda$ defining the equidistant to vary, the generic singularities
are now $A_1$ (smooth surface), $A_2$ (cusp edge), $A_3$ (swallowtail), $A_3^\pm$ (swallowtail beaks/lips transition),
$A_4$ (butterfly) and also $D_4^\pm$ (purse/pyramid) (compare \cite{GZ0}).
The context of the present paper is to extend this to 1-parameter families of surfaces, the parameter in the family
being $\varepsilon$ in our notation, so that there are now two parameters to consider, $\lambda$ and $\varepsilon$.
The particular degeneracy in the $\varepsilon$ family studied here comes from a
`supercaustic chord', that is a chord joining two parabolic points with parallel tangent planes and parallel asymptotic directions.
This occurs
generically only in a 1-parameter family of surfaces.  Along such a chord there may be special values of $\lambda$ where
singularities become more  degenerate, depending on the relative local geometry of the surface patches at the ends of the chord.
When two such special values exist (our Case~1.2) this  corresponds to the intersection of an $A_3$ stratum with the supercaustic.
In addition, there always exists a value of $\lambda$, which we call the degenerate Case~2. This corresponds to the
intersection of a $D_4$ stratum with the supercaustic, and we elucidate ten geometrically distinct cases.
Our paper also gives a natural geometric setting
for  many singularity types which belong to the list of corank~1 maps from $\RR^3$ to $\RR^3$ (\cite{Marar-Tari,Gory}), with the addition of a quadratic term in the extra variable
which does not affect the critical set.  The cases where equidistants are defined by $\lambda =0 $ or 1 remain to be studied.

A second natural 1-parameter family of surfaces is derived from  the `tangential' case in which two surface pieces share a common tangent plane (see for example \cite{GZ});
here boundary singularities occur in the generic case, so that  making one contact point parabolic in a 1-parameter family will introduce additional boundary singularities.
The full  adjacency diagram for singularities of equidistants of 1-parameter families of surfaces, not restricted to the supercaustic case, also remains to be found.

\bigskip
\noindent
{\sc Acknowledgement} \
We are grateful to Aleksandr Pukhlikov for helpful discussions on calculating self-intersections.

\noindent
Peter Giblin, Department of Mathematical Sciences, The University of Liverpool, Liverpool L69 7ZL, UK, email pjgiblin@liv.ac.uk\\
Graham Reeve, Department of Mathematics and Computer Science, Liverpool Hope University, Liverpol L16 9JD, UK, email reeveg@hope.ac.uk

\end{document}